\documentclass[letterpaper,10pt,english]{smfart}
\usepackage{amsmath,amssymb,amsthm,amsxtra,amscd,amsfonts,amstext,latexsym}
\usepackage[T1]{fontenc}
\usepackage{graphicx}
\usepackage{subfigure}
\usepackage[all]{xy}
\usepackage{mathpazo}
\usepackage{hyperref}
\usepackage[lmargin=2.97cm, rmargin=2.97cm, tmargin=2.75cm, bmargin=2.75cm]{geometry}
\newtheorem{theorem}{Theorem}[section]
\newtheorem{lemma}[theorem]{Lemma}
\newtheorem{proposition}[theorem]{Proposition}
\newtheorem{corollary}[theorem]{Corollary}
\theoremstyle{definition}
\newtheorem{remark}[theorem]{Remark}
\newtheorem{example}[theorem]{Example}
\newtheorem{definition}[theorem]{Definition}
\title[Finite symmetry groups in complex geometry]
{Finite symmetry groups in complex geometry}
\author{Kristina Frantzen and Alan Huckleberry%
\thanks{
First
    author's research supported by grants from the Studienstiftung des
    deutschen Volkes and the Deutsche Forschungsgemeinschaft.\\
Second author's research partially supported by grants from the Deutsche
    Forschungsgemeinschaft.}%
}
\address{Institut und Fakult\"at f\"ur Mathematik\\
Ruhr-Universit\"at Bochum\\
Germany\\
kristina.frantzen@ruhr-uni-bochum.de\\
ahuck@cplx.ruhr-uni-bochum.de}
\date { }
\begin{document}
\renewcommand{\baselinestretch}{1.1}
\maketitle
%
%
\subsection* {Introduction}
On June 5, 2007 the second author delivered a talk at
the Journ\'ees de l'Insti\-tut \'{E}lie Cartan entitled 
\emph{Finite symmetry groups in complex geometry}. This
paper begins with an expanded version of that talk
which, in the spirit of the Journ\'ees, is intended
for a wide audience. The later paragraphs are devoted
both to the exposition of basic methods, in particular
an equivariant minimal model program for surfaces,
as well as an outline of recent work of the authors
on the classification of K3-surfaces with special
symmetry.
%
%
\section {Riemann surfaces}
Throughout this introductory section $X$ denotes a Riemann
surface, namely a connected, compact complex manifold of
dimension one.  If it is regarded as a real surface,
then its
\emph{genus} $g=g(X)$, which can be defined 
by $2g$ being the rank of the first homology group
$H_1(X,\mathbb Z)$, is its only topological invariant. 
This means that two such surfaces are homeomorphic
if and only if they are diffeomorphic and this holds if
and only if they have the same genus. Starting with
the case $g=0$, we consider Riemann surfaces from the
point of view of biholomorphic symmetry, i.e., we are
interested in actions on $X$ of subgroups of the group
$\mathrm {Aut}(X)$ of its holomorphic automorphisms.
\subsection {The Riemann sphere}
If $g(X)=0$, then it can be shown that $X$ is 
biholomorphically equivalent to the Riemann sphere.
As the name indicates, this is the 2-dimensional sphere which
often is regarded as the compactification of the complex plane
by adding the point at infinity.  If $z$ is the standard linear
coordinate in a neighborhood of $0$ in the complex plane, then 
$\zeta :=\frac{1}{z}$ is a coordinate of the corresponding 
neighborhood of $\infty $.  

Here we regard the Riemann sphere as the space 
$\mathbb P_1(\mathbb C)$ of 
1-dimensional linear subspaces (lines) of $\mathbb C^2$.
Every point $(z_0,z_1)\in \mathbb C^2\setminus \{(0,0)\}$
is contained in a unique such line which is denoted
by $[z_0:z_1]$.  It follows that $[z_0:z_1]=[w_0:w_1]$
if and only if there exists $\lambda \in \mathbb C^*$
with $(w_0,w_1)=\lambda (z_0,z_1)$. Thus $\mathbb P_1 = \mathbb P_1(\mathbb C)$
can be regarded as the quotient of $\mathbb C^2\setminus \{(0,0)\}$
by scalar multiplication equipped with
the quotient topology. The sets $U_i:=\{[z_0:z_1] \in \mathbb P_1 \, |\,  z_i\not=0\}$
are open. On $U_0$, $U_1$ respectively, we define the coordinate
$z$ by $[z_0:z_1]\mapsto \frac{z_1}{z_0}=:z$,
respectively $[z_0:z_1]\mapsto \frac{z_0}{z_1}=:\zeta $.
Setting $0=[1:0]$ and $\infty =[0:1]$, the description of the Riemann sphere as the set of lines in $\mathbb C^2$ is seen to coincide with the compactified complex plane.

A linear map $T\in \mathrm {GL}_2(\mathbb C)$ takes lines to
lines and therefore induces a map 
$\mathbb P(T):\mathbb P_1\to \mathbb P_1$.  One checks that
$\mathbb P(T)$ is holomorphic and that the homomorphism
$\mathbb P:\mathrm {GL}_2(\mathbb C)\to \mathrm {Aut}(X)$
is surjective.  The kernel of $\mathbb P$ is the
group $\mathbb C^*=\mathbb C^*\cdot \mathrm {Id}$.  Consequently,
if we restrict $\mathbb P$ to $\mathrm {SL}_2(\mathbb C)$,
it is still surjective and yields the exact sequence
$$
0\to \{\pm \mathrm {Id}\}\to \mathrm {SL}_2(\mathbb C)
\to \mathrm {Aut}(\mathbb P_1)\to 0\,.
$$
The subgroup $\{\pm \mathrm {Id}\}$ is the center
of $\mathrm {SL}_2(\mathbb C)$ and
the quotient 
$\mathrm {SL}_2(\mathbb C)/\{\pm \mathrm {Id}\}
=:\mathrm {PSL}_2(\mathbb C)
\cong \mathrm {Aut}(\mathbb P_1)$
is the associated \emph{projective linear group}.

As a subgroup of 
$\mathrm {SL}_2(\mathbb C)$
the special unitary group $\mathrm {SU}_2$
defined by the standard Hermitian structure on $\mathbb C^2$
acts on $\mathbb P_1$.  Since it contains the
center $\{\pm \mathrm {Id}\}$ which acts trivially, this defines
an action of the quotient $\mathrm {SU}_2/\{\pm  \mathrm{Id}\}$,
which can be identified with the group $\mathrm {SO}_3(\mathbb R)$
of orientation preserving linear isometries of $\mathbb R^3$.

Guided by our interest in finite symmetry groups we consider 
finite subgroups of $\mathrm {Aut}(\mathbb P_1 )$.
For the moment we simplify the discussion and
only consider finite subgroups of $\mathrm {SL}_2(\mathbb C)$.
Note that if $G$ is such a subgroup, then we may average
the standard Hermitian structure to obtain a $G$-invariant
Hermitian form and consequently $G$ is conjugate in
$\mathrm {SL}_2(\mathbb C)$ to a subgroup of $\mathrm {SU}_2$.  
If we perform this conjugation,
which changes nothing essential, and project $G$
to $\mathrm {Aut}(\mathbb P_1)$, we may regard
it as a group of Euclidean isometries of $S^2$.  Conversely,
paying the price of the 2:1 central extension, we may
consider the preimage of a group of Euclidean motions in $\mathrm {SU}_2$ and
regard it as acting by holomorphic transformations on
$\mathbb P_1$.

If $X$ is a Riemann surface and $G\subset \mathrm {Aut}(X)$
is a finite group, then the quotient $X/G$ carries a unique
structure of a Riemann surface with the property that
the quotient map $X\to X/G$ is holomorphic.  Following ideas of
Felix Klein, we begin with
a finite group $G$ of rigid motions of the sphere, lift
it to a group of holomorphic transformations of
$X=\mathbb P_1$ in $\mathrm {SU}_2$ and consider such a quotient.
Since there is no nonconstant holomorphic map
from $\mathbb P_1$ to some other Riemann surface, it follows
that $X/G$ is likewise $\mathbb P_1$. Using this fact and
looking closely at the ramified covering map $X\to X/G$,
Klein listed \emph{all} possible finite subgroups of $\mathrm {Aut}(\mathbb P_1)$. Other than
the cyclic and dihedral groups, these are the isometry groups
of the tetrahedron, the octahedron and the icosahedron of order 12, 24 and 60.

The cyclic group $C_n$ of order $n$ can be realized
as a group of diagonal matrices (rotations) in 
$\mathrm {SU}_2$. The dihedral group $D_{2n}$ is 
a semidirect product $C_2\ltimes C_n$. If $C_n$
is a group of diagonal matrices in $\mathrm {SU}_2$, then conjugation with 
\begin {gather*}
w=
\begin {pmatrix}
0 & -1\\
1 & 0
\end {pmatrix}
\end {gather*}
acts on $C_n$ by $x\to x^{-1}$ and $\langle w \rangle \ltimes C_n$ is a realization of $D_{2n}$ in $\mathrm {SU}_2/\{\pm \mathrm{Id}\}$. 

It is a rather simple matter
to place the corner points of a regular tetrahedron and of
a regular octahedron on $\mathbb P_1$ and then to write down
the matrices in $\mathrm {SU}_2$ which realize their
isometry groups as subgroups of $\mathrm {PSU}_2$.  The same
can be done for the iscosahedron, but this is a much more
difficult task. The group of isometries of the
icosahedron is isomorphic to the alternating group $A_5$.
It should be emphasized that the
preimage in $\mathrm{SU}_2$ of a group of rigid motions of a regular polyhedron is a \emph{nontrivial} central
extension.  Particularly in the case of $A_5$, it is an interesting exercise to find the 2-dimensional representation 
of this group!
\subsection {Tori}
In order to identify a Riemann surface $X$ of
genus zero with $\mathbb P_1$, one constructs a 
meromorphic function on $X$ which has only one pole and that being of order one. 
Analogously, one would like to identify a Riemann surfaces of genus one with a complex torus. In order to do this, 
one must prove the existence of
a nowhere vanishing holomorphic 1-form, i.e.,
a 1-form which in local coordinates is given by 
$fdz$ for a nowhere vanishing holomorphic
function $f$.  Integrating this form provides a biholomorphic
map $\alpha :X\to \mathbb C/\Gamma $, where the lattice
$\Gamma $ is the additive subgroup of $(\mathbb C,+)$
defined by integrating the given 1-form over the closed
curves in $X$.  

Unlike the case of $\mathbb P_1$, there is a 1-dimensional
family of holomorphically inequivalent Riemann surfaces
of genus one.  One way of realizing this family is
to choose a basis of the periods so that
$\Gamma =\langle 1,\tau \rangle$, where $\tau $ is in
the upper-halfplane $H^+$.  Then the 1-dimensional family
is given by $H^+/\mathrm {SL}_2(\mathbb Z)\cong \mathbb C$.

A \emph{torus} $T=\mathbb C/\Gamma $ is a group
and acts on itself by group multiplication. This defines
an embedding $T\hookrightarrow \mathrm {Aut}(T)$. 
Given any holomorphic automorphism $\varphi \in \mathrm {Aut}(T)$ 
we lift it to a biholomorphic map 
$\hat \varphi :\mathbb C\to \mathbb C$ of the universal
cover of $T$.
Since $\hat\varphi $ is an affine map of the form 
$\hat\varphi (z)=az+b$ and we are free to conjugate
it with a translation, we may assume that
it is either a linear map, i.e., $z\mapsto az$, 
or a translation.  The translations correspond to
the action of $T$ on itself mentioned above.

Every torus possesses the holomorphic automorphism 
$\sigma $ defined by $t\mapsto -t$, which is also a group automorphism. Except for two special tori,
the full automorphism of $T$ is just $T\rtimes \langle \sigma \rangle$.
These special tori are defined by the lattices $\langle 1,i\rangle $ and
$\langle 1,e^{i\frac{\pi}{6}} \rangle$.  In the
former case $\mathrm {Aut}(T)=T\rtimes C_4$, where the linear
part $C_4$ is generated by a rotation by $\frac{\pi}{2}$, and
in the latter case $\mathrm {Aut}(T)=T\rtimes C_6$, where the
$C_6$ is generated by rotation through 60 degrees.

In summary, in all cases $\mathrm {Aut}(T)=T\rtimes L$ for a linear group $L$ of rotations.  Hence, given a finite subgroup
$G\subset \mathrm {Aut}(T)$, we can decompose it into its
translation and linear parts.
\subsection {Riemann surfaces of general type}
For the remaining Riemann surfaces, i.e., for most,
we have the following observation.
\begin {theorem}\label {finiteness}
The automorphism group of a Riemann surface of genus at
least two is finite.
\end {theorem}
To prove this theorem one needs basic results on the existence
of certain globally defined holomorphic tensors.
For example, one knows that
the space $\Omega (X)$ of holomorphic 1-forms is
$g$-dimensional. 
Note that if $X$ posseses a nowhere vanishing holomorphic 1-form $\omega _0$, then 
every other
holomorphic 1-form $\omega $ is a multiple $\omega =f\omega _0$
where $f$ is a globally defined holomorphic function and therefore constant. Thus 
$\Omega (X)=\mathbb C\omega _0$ and $g(X)=1$. Conversely, if $g(X) >1$ then
every holomorphic 1-form vanishes at at least one point of $X$.
It can be shown that, counting multiplicities,
every $\omega \in \Omega (X)$ has exactly $2g(X)-2$
zeros.

Another basic fact which is useful for the proof of the above
theorem is that the group of holomorphic automorphisms of
a compact complex manifold, in this case a Riemann surface, 
is a \emph{complex Lie group} acting holomorphically on
$X$.  This means that $\mathrm {Aut}(X)$ is itself a 
(paracompact) complex manifold having the property that 
the group operations and the action map
$\mathrm {Aut}(X)\times X\to X$, $(g,x)\mapsto gx$,
are holomorphic. Note that if $\mathrm {Aut}(X)$ is
positive-dimensional and $\{g_t\}$ is a holomorphic 
1-parameter subgroup, then differentiation with respect
to this group defines a holomorphic vector field on $X$.
Conversely, holomorphic vector fields on \emph{compact}
complex manifolds can be integrated to define 
1-parameter groups.

If $\omega \in \Omega (X) \backslash \{0\}$ and $\xi $
is a holomorphic vector field which is not identically
zero, then $\omega (\xi )$ is a holomorphic function on
$X$ which is also not identically zero.  Thus, if
$\omega $ vanished at some point of $X$, we
would have produced a nonconstant holomorphic function,
contrary to $X$ being compact. As a result we obtain the following weak
version of the above theorem.
\begin {proposition}
The automorphism group of a Riemann surface of genus
at least two is discrete.
\end {proposition}
\begin{proof}[Proof of Theorem \ref{finiteness}] One can show that
a Riemann surfaces $X$ of genus
at least two possess enough holomorphic forms,
or holomorphic tensors of higher order, locally
of the form $f(dz)^k$ with $f$ holomorphic, to
definine a canonical embedding of $X$ in a projective space:
If $V_k$ is the vector space of such $k$-tensors, then one
considers the holomorphic map
$
\varphi _{k}:X\to \mathbb P(V_k^*)$
which is defined by sending a point $x\in X$ to the 
\emph{hyperplane} $H_x$ of tensors in $V_k$ vanishing
at $x$.  For $k$ large enough  $\varphi_k$ is
a holomorphic embedding. In fact $k=1$ is usually
enough and at most $k=3$ is required.

The image $Z_k:=\varphi _k(X)$ is a complex submanifold
of the projective space $\mathbb P(V_k^*)$.
Applying Chow's theorem, it follows that $Z_k$ is a algebraic submanifold, i.e., it
is defined as the common zero-set of finitely many homogeneous
polynomials.

Since 
$\mathrm {Aut}(X)$ acts as a group of linear 
transformations on $V_k$, where the action is
given by a representation 
$\rho :\mathrm {Aut}(X)\to \mathrm {GL}(V_k^*)$,
it follows that $\varphi _k$ is $\mathrm {Aut}(X)$-equivariant.
In other words, for every $g\in \mathrm {Aut}(X)$, it follows that
$\varphi _k(gx)=\rho (g)(\varphi _k(x))$.  Thus 
$\mathrm {Aut}(X)$ can be regarded as the stabilizer of $Z_k$
in the projective linear group $\mathrm {PGL}(V_k^*)$.  Since
stabilizers of algebraic submanifolds are algebraic groups
and algebraic groups have only finitely many components,
it follows that $\mathrm {Aut}(X)$ is finite.
\end{proof}
\begin{remark}
In complex geometry the terminology ``manifold
of general type'' refers to a compact complex manifold 
(usually algebraic) which has as many holomorphic tensors of a
certain kind as possible.  In the higher-dimensional case one considers
holomorphic volume forms which are locally of the form
$fdz_1\wedge \ldots \wedge dz_n$ and higher-order tensors
which are described in local coordinates by 
$f(dz_1\wedge \ldots \wedge dz_n)^k$.  One cannot quite require
that the space $V_k$ defines an embedding as above, but it does
make sense to require that the analogous map $\varphi _k$ 
is bimeromorphic onto its image.  Such maps are embeddings outside
small sets.  This is the origin of our referring
to Riemann surfaces of genus at least two as being of general type.
\end{remark}
\subsection {The Hurwitz estimate} 
As in the previous section, we restrict our considerations to
Riemann surfaces $X$ of genus at least two. Having shown that 
$\mathrm {Aut}(X)$ is finite we would like to outline some
ideas behind the proof of the following beautiful theorem.
\begin {theorem}
If $X$ is a Riemann surface of general type, then
$$
\vert \mathrm {Aut}(X)\vert \le 84(g-1)\,.
$$
\end {theorem}
Before going into the ideas of the proof, we emphasize the qualitative meaning of this estimate: the topological Euler
number of $X$ is given by $e(X)=2-2g(X)$ and consequently the 
estimate above is given by $-42e(X)$.  In other words,
the bound for 
$\vert \mathrm {Aut}(X)\vert $ is a linear function of the topological Euler number.

The key to the above estimate is the Riemann-Hurwitz formula
which in our particular case of interest gives a precise
relationship between the topological Euler numbers of
$X$ and $X/G$, where $G$ is any finite group of automorphisms.
In the figure below we have shown a possible
example where the group $G$ is the cyclic group $C_6$ of order six. 
The surface $X$ is schematically represented by a collection
of curves which come together at a number of ramification 
points. The map from upstairs to downstairs represents
the quotient $\pi:X\to X/G$. The observation that with three exceptions
the preimage of a point downstairs consists of six different
points reflects the following general fact: If $G$ is a finite
group in $\mathrm {Aut}(X)$, then there is a finite subset
$R$ such that $G$ acts freely on the complement $X\setminus R$.
\begin{figure}[h]
\begin{center}
\includegraphics[width = 0.5\textwidth]{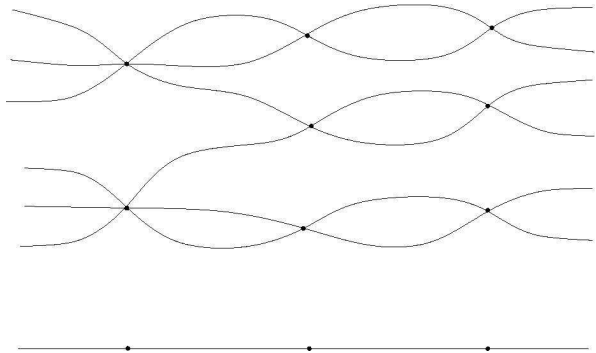}
\caption{A ramified covering of Riemann surfaces}
\end{center}
\end{figure}

At a \emph{ramification point} $x\in R$ the isotropy group
$G_x = \{ g \in G \, | \, g.x=x\}$ consists of more than just the identity.  We note
that the natural representation of $G_x$ on the holomorphic tangent space
$T_xX$ is faithful, and, since 
$\mathrm {GL}(T_xX)\cong \mathbb C^*$, it follows that
$G_x$ is cyclic. Let $\pi (R)=:B$ 
denote the \emph{branch set} of the covering and note
that for every point $b\in B$ we have the canonically defined
numerical invariant $n(b):=\vert G_x\vert $, where $x$ is
any point in the preimage $\pi ^{-1}(b)$. In the figure $B$
consists of three points, two of which have $n(b)=2$ and one
of which has $n(b)=3$.

Let us compare the topological Euler numbers of $X$ and $X/G$
by triangulating $Y:=X/G$ so that set $B$ of branch points is
contained in the set of vertices of the triangles and let us
lift this triangulation to $X$.  We compute the Euler 
number as $v-e+f$, where $v$ is the total number of vertices,
$e$ is the number of edges and $f$ is the number of faces
in the triangulation.  In general,
$$
e(X)=\vert G\vert \cdot e(Y)-\varepsilon _R\,,
$$
where $\varepsilon _R$ is a correction caused by ramification.
In the case of the figure above, every face and every edge
of the triangulation $Y$ lifts to 6 faces and 6 edges 
in the triangulation of $X$.  This is also true for the
vertices which are not contained in $B$.  
In each of the two
cases where $n(b)=2$ we must subtract a correction term
$3=3(2-1)=\vert G.x\vert \cdot  (n(b)-1)$. In the case of $n(b)=3$ 
this correction equals $4 = 2(3-1)$ and the precise formula for
the figure is $e(X)=6e(X/G)-4 - 3- 3$.
In general, the Riemann-Hurwitz formula reads
\begin {align*}
2-2g(X)&=\vert G\vert \cdot(2-2g(Y))  -\sum _{b \in B}\frac{\vert G\vert}{n(b)}(n(b)-1)\\
&=\vert G\vert \cdot  \Bigl((2-2g(Y))-\sum_{b \in B} (1-\frac{1}{n(b)})\Bigr)
\end {align*}
The Hurwitz estimate for the maximal
order of $G$, i.e., for the order of the full automorphism
group $G=\mathrm {Aut}(X)$ follows from experiments with the numbers $n(b)$. (see e.g.\,\cite {K}, Theorem III.2.5)
\subsection {Plane curves}  
One important class of Riemann surfaces consists of those which can
be realized as submanifolds of 2-dimensional projective space.
These are often simply referred to as (complex) \emph{curves}.
Being 1-codimensional, a curve $C$ can be described as the zero-set
of a homogeneous polynomial which is unique up to a scalar factor and vanishes
along $C$ of order one. If $C=\{[z_0:z_1:z_2] \, | \, P_d(z_0,z_1,z_2)=0\}$
where $P_d$ is of degree $d$, then $C$ is said to be of degree
$d$.  Remarkably, the genus of $C$ can be directly computed
from its degree:
$$
g(C)=\frac{(d-1)(d-2)}{2}\,.
$$

\begin{example}
Curves of degree three are of genus one, i.e., they
can be realized as tori $C=\mathbb C/\Gamma $ which have 
positive-dimensional automorphism groups. Very few of these
automorphisms can be realized as restrictions of automorphisms
of $\mathbb P_2$ which, in analogy to the case of $\mathbb P_1$,
are induced by linear transformations of $\mathbb C^3$.  
To see this, let $S:=\mathrm {Stab}_{\mathrm {Aut}(\mathbb P_2)}(C)$
be the subgroup of elements $T\in \mathrm {Aut}(\mathbb P_2)$
with $T(C)=C$.  Since $C$ is defined by a complex polynomial
equation, the group $S$ is a complex subgroup of the complex Lie group
$\mathrm {Aut}(\mathbb P_2)=\mathrm {SL}_3(\mathbb C)/C_3$. 
Here $C_3$ is realized in $\mathrm {SL}_3$ as its center,
i.e, the group of diagonal matrices 
$ \lambda \cdot  \mathrm{Id}_{\mathbb C^3}$ with $\lambda ^3=1$.

Note that no element of $S$ fixes $C$ pointwise, because
such a linear transformation would necessarily pointwise fix
the linear subspace of $\mathbb P_2$ spanned by $C$, i.e.,
$\mathbb P_2$ itself.  Since $\mathrm {Aut}(C)$ is compact,
it follows that $S$ is likewise compact. Thus the lift
$\hat S$ of $S$ to $\mathrm {SL}_3(\mathbb C)$ can be
regarded as a compact complex submanifold of the vector space
$\mathrm {Mat}(2\times 2,\mathbb C)\cong \mathbb C^4$. Consequently,
we obtain holomorphic functions on $\hat S$ as restrictions
of holomorphic functions on $\mathbb C^4$. Since $\hat S$ 
is compact, the maximum
principle implies that these are
constant on its components. As these restrictions clearly
separate the points of $\hat S$ we
see that $\hat S$ is finite. Hence we have proven the following proposition.
\begin {proposition}
If $C$ is a cubic curve in $\mathbb P_2$, then the subgroup
$S$ of $\mathrm {Aut}(C)$ of automorphisms which extend to
automorphisms of $\mathbb P_2$ is finite.
\end {proposition}
\end{example}

The case of a curve $C$ of degree four is completely different.
In this case $g(C)=3$ and the space
$V_1$ of holomorphic 1-forms on $C$ is 3-dimensional. The mapping $\varphi _1:C\to \mathbb P(V_1^*)$ is an embedding
and the orginal realization of $C$ as a curve is, after a choice
of coordinates, just $\mathrm {Im}(\varphi _1)$.  Consequently,
curves of degree four are equivariantly embedded.
\begin {proposition}
If $C$ is a quartic curve, then every automorphism
of $C$ extends to a unique automorphism of $\mathbb P_2$.
\end {proposition}
The study curves of genus
three from the point of view of symmetry 
is therefore closely related to the classification and invariant
theory of finite subgroups of $\mathrm {SL}_3(\mathbb C)$
(\cite {B}, see also \cite {Y})
\begin{example}
We consider the quartic curve
defined by 
$z_0z_1^3+z_1z_2^3+z_2z_0^3$ and refer to it as Klein's curve $C_{\mathrm {Klein}}$. 
Although the general theory tells us
that every automorphism of $C_{\mathrm {Klein}}$ is the restriction
of an automorphism of $\mathbb P_2$, not all of these automorphisms
are immediately visible. In fact,
$
\mathrm {Aut}(C)\cong \mathrm {PSL}_2(\mathbb F_7)
\cong \mathrm {GL}_3(\mathbb F_2)$.
This group, which is often denoted by $L_2(7)$, is the
unique simple group of order 168. 
Note that $168=84(3-1)$
and therefore the automorphism group of Klein's curve
attains the maximal order among Riemann surfaces of genus three allowed by the Hurwitz estimate. One
can show that $C_{\mathrm {Klein}}$ is the unique genus three curve for which this
upper bound is attained. 
The book \cite{8fold} is dedicated to various interesting aspects concerning the geometry of this curve and its automorphisms.
\end{example}
\section {Manifolds of general type}
The title of this work indicates our interest in the role
of finite symmetry groups in arbitrary dimensions.  Nevertheless,
after the previous introductory section on Riemann surfaces,
mostly all of our considerations are devoted to 
the case of compact complex surfaces, i.e., complex 2-dimensional, 
connected, compact complex manifolds.  Before restricting
to that case, we do comment on higher-dimensional manifolds 
of general type.

Recall that a Riemann surface $X$ is of general type if it
possesses sufficiently many globally defined holomorphic tensors,
which are locally of the form $f(dz)^k$, so that the $k$-canonical map
$\varphi _k:X \to \mathbb P(V_k^*)$ is a biholmorphic
embedding for $k$ sufficiently large. In the higher-dimensional
case, $\mathrm {dim}_{\mathbb C}X=n$, the analogous objects in
the case $k=1$ are holomorphic $n$-forms which are
locally of the form $fdz_1\wedge \ldots \wedge dz_n$. Here
$f$ is a holomorphic function on a coordinate chart with
coordinates $z_1,\ldots ,z_n$.  For arbitrary $k$, its $k$-tensors 
are locally of the form
$f(dz_1\wedge \ldots \wedge dz_n)^k$. In a less archaic language,
these are sections of the $k$-th power $\mathcal K_X^k$ of the canonical line
bundle and the space $V_k$ of $k$-tensors is denoted
by $\Gamma (X, \mathcal K_X^k)$.

It turns out to be appropriate to require that for some $k$
the mapping 
$\varphi _k:X \to \mathbb P(V_k^*)$ is a meromorphic 
instead of holomorphic embedding.
This condition is 
more conveniently described by an invariant of the
associated function field:
if $s$ and $t$ are two tensors of the same type, i.e.,
$s,t\in V_k$, where locally $s=f(dz_1\wedge \ldots \wedge dz_n)^k$
and $t=g(dz_1\wedge \ldots \wedge dz_n)^k$, then their
ratio 
$$
m=\frac{s}{t}=\frac{f(dz_1\wedge \ldots \wedge dz_n)^k}
{g(dz_1\wedge \ldots \wedge dz_n)^k}=\frac{f}{g}
$$
has an interpretation as a globally defined
meromorphic function on $X$.  

We let $Q(V_k)$ be the
quotient field generated by this procedure. From the
theorem of Thimm-Siegel-Remmert we know that
meromorphic functions on a compact complex manifold are analytically
dependent if and only if they are algebraically dependent, and consequently the transcendence
degree of $Q(V_k)$ over the field of constant functions, or equivalently, the maximal number of analytically independent
meromorphic functions which can be constructed as quotients
of tensors from $V_k$, 
is at most $\mathrm {dim}_{\mathbb C}(X)$.
If for some $k$ this number equals
the dimension of $X$, then one says that $X$ is of 
\emph{general type}.  

If $X$ is of general type, then
for some $k$, maybe not the one in the definition, the
mapping $\varphi _k$ is indeed a bimeromorphic 
embedding.  Since $\varphi _k$ is equivariant with respect
to the full automorphism group, it follows that 
$\mathrm {Aut}(X)$ is represented as the subgroup of
$\mathrm {Aut}(\mathbb P(V_k^*))$ which stabilizes the
image $\mathrm {Im}(\varphi _k)$.  Thus, by precisely
the same type of argument as in the 1-dimensional case,
we have the following fact.
\begin {proposition}
The automorphism group of a manifold of general type
is finite.
\end {proposition}
Even if the field $Q(V_k)$ does not have transcendence
degree equal to the dimension of $X$, it
certainly contains important information. As a very
rough first invariant one defines the \emph{Kodaira dimension}
of $X$ as the maximal transcendence degree attained by 
some $V_k$.  In other words, the Kodaira dimension $\kappa (X)$
is the maximal number of analytically independent meromorphic
functions which can be obtained as quotients of $k$-tensors.
If for all $k$ there are no such tensors, i.e., $V_k=0$
for all $k$, then one lets $\kappa (X):=-\infty$.  It follows
that $\kappa (X)\in \{-\infty, 0,1,\ldots ,\mathrm {dim}(X)\}$.

For 
$\mathrm{dim}(X) >\kappa (X)\ge 1$ the meromorphic maps $\varphi _k:X\to \mathbb P(V_k^*)$
can still be very interesting.
However, it is quite possible
that a nontrivial automorphism 
of $X$ can act trivially on $\mathrm {Im}(\varphi _k)$.
\subsection {Surfaces of general type - Quotients by small subgroups}\label{smallgroups}
Inspired by the 1-dimen\-sional case one wishes to obtain bounds for the order of the automorphism groups of manifolds of general type.
Below we give an outline of a simple method which has been used, for example, to
obtain estimates of Hurwitz type (see \cite{HS}). 

At the 
present time the following
sharp estimate for surfaces of general type requires essentially
more combinatorial work (\cite {XG}, \cite{XGa}).  
\begin {theorem}
The order of the automorphism group of a surface of
general type is bounded by $(42)^2K^2$.
\end {theorem}
Since the self-intersection number $K^2$ of a canonical divisor is a topological
invariant, this is exactly the desired type of estimate.

Turning to the method mentioned above, 
given a finite group $G$ we
want analyze the possibilities of it acting on surfaces
of general type with given topological invariants, in this
case Chern numbers. We then look for a small subgroup
$S$ in $G$ with an interesting normalizer $N$. The
notion of interesting can vary. For example, this can mean that $N$ is large with respect to $G$ or that $N$ has a rich group structure. The group $S$ should
be small in the sense of size and structure. Clearly,
$S=C_2$ or some other small cyclic group would be 
a good choice.

The first step is to consider the quotient $X\to X/S=:Y$.
Since the normalizer $N$ acts on $Y$, we are presented
with the new task of understanding $Y$ as an $N$-variety
and $X\to Y$ as an $N$-equivariant map, e.g.\,by studying the action of $N$
on the ramification and branch loci. If this can be
done, then we attempt to piece together the $G$-action
on $X$ from knowledge of the $S$-quotient and the
$N$-action on $Y$. 

Because we have
been forced to transfer our consideration to the smaller group $N$, it might appear that we have even lost ground. However,
there are at least two possible advantages of this approach.
First, without being overly optimistic one can hope that
the topological invariants have decreased in size so that
if $Y$ is still of general type, some inductive argument
can be carried out. Alternatively, if $Y$ is not of general type,
then we come into a range of Kodaira dimension where new methods are available.  
%
%
%
%
\section {The Enriques Kodaira classification}\label{EKC}
From now on we restrict our considerations to compact complex
surfaces $X$. Here $\kappa (X)$ can take on the values $-\infty $,
$0$, $1$ and $2$. One would like to prove a classification theorem
similar to that for Riemann surfaces with one big class
consisting of the surfaces of general type and the remaining
surfaces with $\kappa (X) \le 1$ being precisely described.

This is almost possible with the final result
being called the \emph{Enriques Kodaira classification}. For a detailed exposition
we refer the reader e.g.\,to \cite {BPV}.  
In the
case of algebraic surfaces, i.e., those compact complex
surfaces which can be holomorphically embedded in some
projective space, much of the essential work was carried
out by members of the Italian school of algebraic geometry,
in particular by Enriques. It should be noted that a
surface is algebraic if and only if it possesses two 
analytically independent meromorphic functions.
\subsection* {Minimal models}
One complicating factor in the classification theory
is that, given a surface $X$, one can \emph{blow it up}
to obtain a new surface $\hat X$ and a holomorphic
map $\hat X\to X$ which is almost biholomorphic.  
Conversely, given $X$, one may be able to \emph{blow it down}, i.e., $X$ is the blow up of some other surface.
Let us briefly explain this process.

The simplest \emph{blow up} is constructed as follows.
Let $X=\mathbb P_2$ and choose $p:=[1:0:0]$ as the point
to be blown up. The projection $\pi :\mathbb P_2\to \mathbb P_1$,
$[z_0:z_1:z_2]\mapsto [z_1:z_2]$ is well-defined
and holomorphic outside of the base point $p$. So we
remove $p$ and consider the restricted map.  Its fiber
over a point $[a,b]$ is a copy of the complex plane parameterized by $t\mapsto [t:a:b]$. 
This continues
to a map of the full projective line given by
$[t_0:t_1]\mapsto [t_0:t_1a:t_1b]$.  
The map $\pi$ realizes $\mathbb P_2\setminus \{p\}$ as a bundle of lines: we say that the fibration 
$\mathbb P_2\setminus \{p\}\to \mathbb P_1$ is a
\emph{line bundle}, a rank one
holomorphic vector bundle. Each of its fibers is naturally isomorphic
to $\mathbb C$ and contains a natural choice of zero,
namely $[0:a:b]$, and a natural point at infinity $[1:0:0]$.

The point at infinity is the same 
for each of the lines.  To resolve this problem, we 
formally add individual points at infinity to each of the lines, i.e., two different lines receive different points at
infinity.  One checks that this construction results in a complex manifold
$\mathrm {Bl}_p(\mathbb P_2)$ to which the line bundle fibration
extends as a $\mathbb P_1$-bundle 
$\mathrm {Bl}_p(\mathbb P_2)\to \mathbb P_1$.  The set $E$ of points
at infinity is a copy of $\mathbb P_1$ which is mapped to
the point $p$ by the natural projection 
$\mathrm{Bl}_p(\mathbb P_2)\to \mathbb P_2$.  Outside of 
$E$ this projection is biholomorphic.

There are several first observations about this construction.
For one, it should be noted that the construction is local.
In other words, for any open neighborhood $U$ of $p$ we can define the blow up $\mathrm {Bl}_p(U)\to U$.  Thus, for
any surface $X$ and a point $p\in X$ we have 
$\mathrm {Bl}_p(X)\to X$.  One checks that up
to biholomorphic transformations the construction is
independent of the coordinate chart which is used.
Secondly, regarding $E$ as a homology class in
$H_2(\mathrm {Bl}_p(X),\mathbb Z)$ which is equipped
with its natural intersection pairing, one shows that
$E\cdot E=-1$.  Remarkably, the converse statement holds
(see \cite{G}):
\begin {theorem}
Let $X$ be a complex surface and $E$ be a smooth curve
in $X$ which is holomorphically equivalent to $\mathbb P_1$.
Then there is a holomorphic map $X\to Y$ to a complex surface
which realizes $X$ as $\mathrm {Bl}_p(Y)$ if and only
if $E\cdot E=-1$.
\end {theorem}
As a result of this theorem it is reasonable to classify only those
surfaces which are \emph{minimal} in the sense that they
contain no (-1)-curves, i.e., curves $E$ which are biholomorphic
to $\mathbb P_1$ and satisfy $E\cdot E=-1$. We give a rough
summary of this classification: Any minimal surfaces belongs to one of the following classes of surfaces, ordered according to Kodaira dimension.

\smallskip
\begin{center}
\begin{tabular}{ll}
 $\kappa = -\infty$ & Ruled surfaces, $\mathbb P_2$, and exceptional nonalgebraic
surfaces\\
$\kappa = 0$ & Tori,  K3-, Enriques-, Kodaira-, and bi-elliptic surfaces\\
$\kappa = 1$ & Elliptic surfaces\\
$\kappa = 2$ & Surfaces of general type
\end{tabular} \end{center}\smallskip
A ruled surface $X \neq \mathbb P_1 \times \mathbb P_1$ admits a canonical locally trivial holomorphic fibration onto a Riemann surface with generic
 fiber $\mathbb P_1$. Elliptic surfaces possess canonically
 defined fibrations over $\mathbb P_1$ with the generic fiber
 being a 1-dimensional torus. In this situation different fibers can be
 biholomorphically different tori.  
Note that canonically 
defined fibrations $\pi:X\to Y$ are automatically equivariant,
 i.e., there is an action of $\mathrm{Aut}(X)$
 on $Y$ so that for every $g\in \mathrm{Aut}(X)$ it follows
 that $\pi \circ g=g\circ \pi $. Thus from the point of view of group actions
 it is particularly advantageous if the
surface is either ruled or elliptic. To exemplify this, we present an detailled discussion of automorphisms of rational ruled surfaces, the Hirzebruch surfaces. 
\subsection*{Hirzebruch surfaces}
The $n$-th Hirzebruch surface $\Sigma _n$ is defined as the compactification
of the total space of the $n$-th power $H^n$ of the hyperplane
bundle over $\mathbb P_1$. The compactification is constructed by adding the
point at infinity to each fiber. This makes sense because
the structure group of a line bundle is 
$\mathrm {GL}_1(\mathbb C)\cong \mathbb C^*$ whose action on
the complex line canonically extends to an action on
$\mathbb P_1$.  The surface $\Sigma _0$ is the compactification
of the trivial bundle and is therefore 
$\mathbb P_1\times \mathbb P_1$. We have seen above that
$\Sigma _1$ is $\mathrm {Bl}_p(\mathbb P_2)$.

By construction the $\mathbb P_1$-bundle
$\Sigma _n\to \mathbb P_1$ has a section $E_n$ at infinity.   
Let us show that $E_n$ can be blown down to a point which,
except in the case of $\Sigma _1$, is singular. For this
it is convenient to recall that $H^n$ is the quotient of
$H$ by the cyclic group $C_n$ acting via the principal
$\mathbb C^*$-action in the fibers of $H$. This extends
to $\Sigma _1$ to give us a diagram
\begin {gather*}
\begin {xymatrix}{
\Sigma _1 \ar[r] \ar[d]^{b_1} & \Sigma _n \ar[d]^{b_n}\\ 
\mathbb P_2 \ar[r]  & \mathrm{Cone}_n(\mathbb P_1).}
\end {xymatrix}
\end {gather*}
The map $b_1$ blows down the (-1)-curve $E_1$. The horizontal maps are the
$C_n$-quotients and $b_n$ is induced by $b_1$. Regarding $\mathbb P_2$ as a cone over the line
$\{z_0=0\}$ at infinity, $\mathrm {Cone}_n(\mathbb P_1)$
is defined as its quotient by $C_n$ acting on its fibers.

Let us turn now to the automorphism groups of the Hirzebruch
surfaces. For this we consider 
the standard $\mathrm {GL}_2(\mathbb C)$-action on $\mathbb P_2$
fixing the point $p = [1:0:0]$ and stabilizing
the hyperplane $\{z_0=0\}$. By construction, it lifts to the
blow up $\Sigma _1=\mathrm {Bl}_{p}(\mathbb P_2)$. This
action is centralized by the $C_n$-action discussed above
and therefore there is an holomorphic action of 
$L_n:=\mathrm {GL}_2(\mathbb C)/C_n$ on the quotient 
$\Sigma _n=\Sigma _1/C_n$.  

The remaining automorphisms come from the sections 
$s \in \Gamma (\mathbb P_1,H^n)$ of $H^n$. In general
if $s:X\to L$ is a section of a holomorphic line bundle $\pi:L\to X$, then 
$x\mapsto x+s(\pi (x))$ defines a holomorphic automorphism
of the bundle space $X$ which extends to the associated 
$\mathbb P_1$-bundle.  Thus, in the case at hand we may
regard $\Gamma (\mathbb P_1,H^n)$ as a subgroup of
$\mathrm {Aut}(\Sigma _n)$. Conjugation by elements
of $L_n$ stabilizes $\Gamma (\mathbb P_1,H^n)$ and 
the semidirect product $L_n\ltimes \Gamma (\mathbb P_1,H^n)$ 
is a subgroup of $\mathrm {Aut}(\Sigma _n)$. In fact there are no other automorphisms
(see e.g.\,\cite {HO}).
\begin {proposition}
$\mathrm {Aut}(\Sigma _n)=L_n\ltimes \Gamma (\mathbb P_1,H^n)$.
\end {proposition}
Having identified $\mathrm {Aut}(\Sigma _n)$ we wish to pin down
its finite subgroups.  First note that the maximal compact subgroups of a connected Lie
group are unique up to conjugation. In the case of
$\mathrm {Aut}(\Sigma _n)$ we observe that the image
of the unitary group $\mathrm {U}_2$ is a maximal compact subgroup.
If $G$ is a finite subgroup of $\mathrm {Aut}(\Sigma _n)$,
then we may conjugate it to a subgroup
of the image of the unitary group. Allowing the action
to contain the kernel of 
$\mathrm {U}_2(\mathbb C)\to \mathrm {U}_2/C_n$, we chararacterize
the finite group actions on $\Sigma _n$ as being given
by finite subgroups of $\mathrm {U_2}$.

\medskip
Let us return the rough classification of minimal surfaces outlined above. 
Note that 
bi-elliptic surfaces admit a locally trivial elliptic fibration over an elliptic curve. The same holds for a 
Kodaira surface itself or an unramified covering.
In analogy to our treatment of Hirzebruch surfaces, automorphisms of fibered surfaces can be investigated by using the structure given by the fibration. 

Finite subgroups acting on $\mathbb P_2$ were classified
at the turn of the 20th century by Blichfeld (\cite {B}, see also \cite{Y}).
Much later the finite groups acting on 2-dimensional tori
were classified by Fujiki (\cite {F}). Noting that the universal cover of an Enriques surface is a K3-surface, the interest of finite symmetry groups may be focussed on the remaining
case of K3-surfaces. 

From the point of view of symmetries the restriction of our 
considerations to minimal surfaces is not necessarily natural. 
This is due to the fact that the reduction from a surface to its minimal model may not be equivariant. The concept of minimal models therefore needs to be replaced by an equivariant analogue, namely an equivariant reduction procedure. The following section is dedicated to a detailed presentation of the 
equivariant minimal model program for surfaces 
formulated in the language of Mori theory. 
%
%
\section {Equivariant Mori reduction}\label{mmp}
This section is dedicated to a discussion of Example 2.18 in \cite{kollarmori} (cf.\ also Section 2.3 in \cite{Mori}) which introduces a minimal model program for surfaces
respecting finite groups of symmetries.

Given a projective algebraic surface $X$ with an action of a finite group $G$, in analogy to the usual minimal model program, one obtains from $X$ a $G$-minimal model $X_{G\text{-min}}$ by a finite number of $G$-equivariant blow-downs, each contracting a finite number of disjoint (-1)-curves. The surface $X_{G\text{-min}}$ is either a conic bundle over a smooth curve, a Del Pezzo surface or has nef canonical bundle. 

Equivariant Mori reduction and the theory of $G$-minimal models have applications in various different contexts and can also be generalized to higher dimensions.
Initiated by Bayle and Beauville in \cite{Bayle}, the methods have been employed in the classification of subgroups of the Cremona group $\mathrm{Bir}(\mathbb P_2)$ of the plane for example by Beauville and Blanc \cite{Beauville}, \cite{BeauBlancPrime}, \cite{PhDBlanc}, \cite{Blanc1}, de Fernex \cite{fernex}, Dolgachev and Iskovskikh \cite{DolgIsk}, \cite{DolgIsk2}, and Zhang \cite{ZhangRational}. 
These references also provide certain relevant details regarding Example 2.18 in \cite{kollarmori} and Section 2.3 in \cite{Mori}: e.g.\ the case $G \cong C_2$ is discussed in \cite{Bayle}, the case $G \cong C_p$ for $p$ prime in \cite{fernex}, the case of perfect fields is treated in \cite{DolgIsk2}.

For the convenience of the reader we here give a detailed exposition of the equivariant minimal model program for arbitrary finite groups acting on complex projective surfaces
(see also Chapter 2 in \cite{Fra} for further details).
\subsection {The cone of curves and the cone theorem}
Throughout this section we let $X$ be a smooth projective algebraic surface and let $\mathrm{Pic}(X)$ denote the 
group of isomorphism classes of line bundles on $X$. Here a curve is an irreducible 1-dimensional subvariety.   
\begin{definition}
A \emph{divisor} on $X$ is a formal linear combination of curves $C = \sum a_i C_i$ with $a_i \in \mathbb Z$.
A \emph{1-cycle} \index{1-cycle} on $X$ is a formal linear combination of curves $C = \sum b_i C_i$ with $b_i \in \mathbb R$. A 1-cycle is \emph{effective} if $b_i \geq 0$ for all $i$. 
Extending the pairing $\mathrm{Pic}(X) \times \{\text{divisors}\} \to \mathbb Z$, $(L,D) \mapsto L \cdot D = \deg(L|_D)$ by linearity, we obtain a pairing $\mathrm{Pic}(X) \times \{\text{1-cycles}\} \to \mathbb R$.
Two 1-cycles $C,C'$ are called \emph{numerically equivalent}  if $L\cdot C = L \cdot C'$ for all $L \in \mathrm{Pic}(X)$.  We write $C \equiv C'$. The numerical equivalence class of a 1-cycle $C$ is denoted by $[C]$.
The space of 1-cycles modulo numerical equivalence is a real vector space denoted by $N_1(X)$. 
Note that $N_1(X)$ is finite-dimensional. A line bundle $L$ is called \emph{nef}\index{nef} if $L \cdot C \geq 0$ for 
all curves $C$. 
We set
\[
NE(X) = \{ \sum a_i[C_i] \ | \ C_i \subset X \text{ irreducible curve},\, 0 \leq a_i \in \mathbb R\} \subset N_1(X).
\]
The closure $\overline{NE}(X)$ of $NE(X)$ in $N_1(X)$ is called \emph{Kleiman-Mori cone} or \emph{cone of curves}. 
For a line bundle $L$, we write $\overline{NE}(X)_{L\geq 0} = \{ [C]\in N_1(X) \ |\ L \cdot C\geq 0 \} \cap \overline{NE}(X)$. Analogously, we define $\overline{NE}(X)_{L\leq 0}$, $\overline{NE}(X)_{L > 0}$, and $\overline{NE}(X)_{L < 0}$.
\end{definition}
Using this notation we phrase Kleiman's ampleness criterion (cf.\,Theorem 1.18 in \cite{kollarmori}) as follows:
A line bundle $L$ on $X$ is ample if and only if 
$\overline{NE}(X)_{L>0} = \overline{NE}(X)\backslash \{0\}$. 
\begin{definition}
A subset $N \subset V$ of a finite-dimensional real vector space $V$ 
is called \emph{cone} if $0\in N$ and $N$ is closed 
under multiplication by positive real numbers. 
A subcone $M \subset N$ is called \emph{extremal} if $u,v \in N $ satisfy $u,v  \in M$ whenever $u+v \in M$. An extremal subcone is also referred to as an \emph{extremal face}. A 1-dimensional extremal face is called \emph{extremal ray}.
For subsets $A, B \subset V$ we define $A+B := \{ a+b \,|\, a\in A, b\in B\}$.
\end{definition}
The cone of curves $\overline{NE}(X)$ is a convex cone in $N_1(X)$ and the following cone theorem, stated here only for surfaces, describes its geometry (cf.\,Theorem 1.24 in \cite{kollarmori}).
\begin{theorem}\label{conethm}
Let $X$ be a smooth projective surface and let $\mathcal{K}_X$ denote the canonical line bundle on $X$. 
There are countably many rational curves $C_i \in X$ such that $0 < -\mathcal{K}_X \cdot C_i \leq \mathrm{dim}(X) +1 $ and 
\[
\overline{NE}(X) = \overline{NE}(X) _{\mathcal{K}_X \geq 0} + \sum_i \mathbb R_{\geq 0} [C_i]. 
\]
For any $\varepsilon>0$ and any ample line bundle $L$
\[
\overline{NE}(X) = \overline{NE}(X) _{(\mathcal{K}_X+\varepsilon L) \geq 0} + \sum_{\text{finite}} \mathbb R_{\geq 0} [C_i]. 
\]
\end{theorem}
\subsection{Surfaces with group action and the cone of invariant curves}
\begin{figure}[h]
\centering
     \subfigure[The cone of curves]
        {\includegraphics[width=0.47\textwidth]{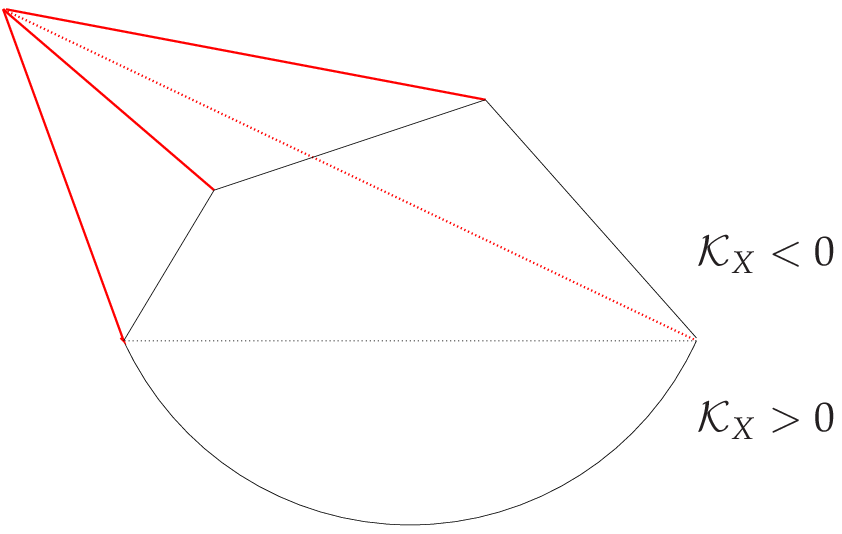}}\hspace{0.8cm}
     \subfigure[The cone of curves and the invariant subspace $N_1(X)^G$. 
                Their intersection $\overline{NE}(X)^G$ has a new extremal ray.]
        {\includegraphics[width=0.47\textwidth]{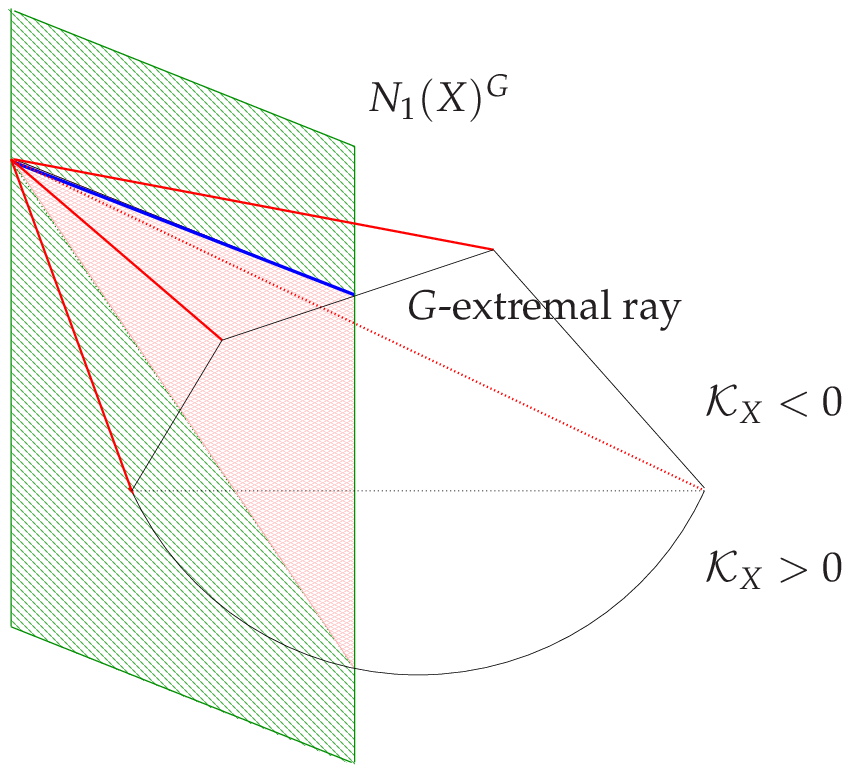}}
      \caption{The extremal rays of $\overline{NE}(X)^G$ are 
                not extremal in $\overline{NE}(X)$}\label{moribild}
\end{figure}

Let $X$ be a smooth projective surface and let $G \subset \mathrm{Aut}(X)$ be a group of holomorphic transformations of $X$. For $g \in G$ and an irreducible curve $C_i$ we denote by $g C_i$ the image of $C_i$ under $g$. For a 1-cycle $C = \sum a_i C_i$ we define $gC = \sum a_i (g C_i)$. This defines a $G$-action on the space of 1-cycles. 
Since two 1-cycles $C_1, C_2$ are numerically equivalent if and only if $gC_1 \equiv gC_2$ for any $g \in G$, we can
define a $G$-action on $N_1(X)$ by setting $g[C] := [gC]$ and extending by linearity. We write $N_1(X)^G = \{ [C] \in N_1(X) \ | \ [C]=[gC] \text{ for all } g \in G\}$, the set of invariant 1-cycles modulo numerical equivalence. This space is a linear subspace of $N_1(X)$.
The cone $NE(X)$ is a $G$-invariant set and so is its closure $\overline{NE}(X)$. The subset of invariant elements in $\overline{NE}(X)$ is denoted by $\overline{NE}(X)^G$.
\begin{remark}
$
\overline{NE}(X)^G = \overline{NE(X) \cap N_1(X)^G}=\overline{NE}(X) \cap N_1(X)^G .
$
\end{remark}
The subcone $\overline{NE}(X)^G$ of $\overline{NE}(X)$ inherits the geometric properties of $\overline{NE}(X)$ established by the cone theorem. 
Note however that the extremal rays of $\overline{NE}(X)^G$, which we refer to as \emph{$G$-extremal rays}, are in general neither extremal in $\overline{NE}(X)$ (cf.\,Figure \ref{moribild}) nor generated by classes of curves but by classes of 1-cycles.

\begin{lemma}\label{Gextremalray} Let $G$ be a finite group and let $R$ be a $G$-extremal ray with $\mathcal{K}_X \cdot R<0$. Then there exists a rational curve $C_0$ such that $R$ is generated by the class of the 1-cycle $C = \sum_{g \in G} gC_0$.
\end{lemma}
\begin{proof}
Consider a $G$-extremal ray $R = \mathbb{R }_{\geq 0}[E]$. By the cone theorem (Theorem \ref{conethm}) $[E] \in \overline{NE}(X)^G \subset \overline{NE}(X)$ can be written as $[E] = [\sum_i a_i C_i] +[F]$,
where $ \mathcal{K}_X \cdot F \geq 0$, $a_i \geq 0$ and $C_i$ are rational curves. Let $|G|$ denote the order of $G$ and let $[GF] = G[F] = \sum_{g\in G} g[F]$. Since $g[E]=[E]$ for all $g\in G$ we can write
\[
|G| [E] = \sum_{g\in G} g[E] = \sum_{g\in G}([\sum_i a_igC_i] + g[F]) = \sum_i a_i G[C_i] + G[F].
\]
The element $[\sum a_i (GC_i)] + [GF]$ of the $G$-extremal ray $\mathbb{R }_{\geq 0}[E]$ is decomposed as the sum of two elements in $\overline{NE}(X)^G$. Since $R$ is extremal in $\overline{NE}(X)^G$ both must lie in $R=\mathbb{R }_{\geq 0}[E]$ .
Consider $[GF] \in R$. Since $g^*\mathcal{K}_X \equiv \mathcal{K}_X$ for all $g \in G$, we obtain
$$\mathcal{K}_X \cdot (GF) = \sum_{g\in G} \mathcal{K}_X \cdot (gF) = \sum_{g \in G}(g^*\mathcal{K}_X) \cdot F= |G|\mathcal{K}_X \cdot F \geq 0.$$ 
As $\mathcal{K}_X \cdot R < 0$ by assumption, this implies $[F]=0$ and $\mathbb{R }_{\geq 0}[E] = \mathbb{R }_{\geq 0}[\sum a_i (GCi)]$. Again using the fact that $R$ is extremal in $\overline{NE}(X)^G$, we conclude that each summand of $[\sum a_i (GC_i)]$ must be contained in $R=\mathbb{R }_{\geq 0}[E]$ and the extremal ray $\mathbb{R }_{\geq 0}[E]$ is therefore generated by $[GC_i]$ for some $C_i$ chosen such that $[GC_i] \neq 0$. This completes the proof of the lemma. 
\end{proof}
\subsection{The contraction theorem and minimal models of surfaces}
In this subsection, we state the contraction theorem for smooth
projective surfaces. The proof of this theorem can be found e.g.\,in
\cite{kollarmori} and needs to be modified slightly in order to give an equivariant contraction theorem in the next subsection. 
\begin{definition}
Let $X$ be a smooth projective surface and let $F \subset \overline{NE}(X)$ be an extremal face. A morphism $\mathrm{cont}_F: X \to Z$ is called the \emph{contraction of $F$} if
\begin{itemize}
\item[-] $(\mathrm{cont}_F)_*\mathcal{O}_X = \mathcal{O}_Z$ and 
\item[-] $\mathrm{cont}_F(C) = \{\text{point}\}$ for an irreducible curve $C\subset X$ if and only if $[C] \in F$. 
\end{itemize}
\end{definition}
The following result is known as the contraction theorem (cf.\,Theorem 1.28 in \cite{kollarmori}).
\begin{theorem}\label{contractionthm} \index{contraction theorem}
Let $X$ be a smooth projective surface and $R \subset \overline{NE}(X)$ an extremal ray such that $\mathcal{K}_X \cdot R<0$. Then the contraction morphism $\mathrm{cont}_R: X \to Z$ exists and is one of the following types:
\begin{enumerate}
\item{$Z$ is a smooth surface and $X$ is obtained from $Z$ by blowing up a point. }
\item{$Z$ is a smooth curve and $\mathrm{cont}_R:X \to Z $ is a minimal ruled surface over $Z$.}
\item{$Z$ is a point and $\mathcal{K}_X^{-1}$ is ample. }
\end{enumerate}
\end{theorem}
The contraction theorem leads to a minimal model program for surfaces: Starting from $X$, if $\mathcal{K}_X$ is not nef, i.e, there exists an irreducible curve $C$ such that $\mathcal{K}_X \cdot C < 0$, then $\overline{NE}(X)_{\mathcal{K}_X <0}$ is nonempty and there exists an extremal ray $R$ which can be contracted. The contraction morphisms either gives a new surface $Z$ (in case 1) or provides a structure theorem for $X$ which is then either a minimal ruled surface over a smooth curve (case 2) or isomorphic to $\mathbb P^2$ (case 3). Note that the contraction theorem as stated above only implies $\mathcal{K}_X^{-1}$ ample in case 3. It can be shown that $X$ is in fact $\mathbb P^2$. This is omitted here since this statement does not transfer to the equivariant setup. In case 1, we can repeat the procedure if $K_Z$ is not nef. Since the Picard number drops with each blow down, this process terminates after a finite number of steps. The surface obtained from $X$ at the end of this program is called a \emph{minimal model} of $X$. 
\begin{remark}
 Let $E$ be a (-1)-curve on $X$ and $C$ be any irreducible curve on $X$. Then $E \cdot C < 0$ if and only if $C =E$. It follows that $\overline{NE}(X) = \mathrm{span}(\mathbb R_{\geq 0}[E], \overline{NE}(X)_{E \geq 0})$. Now $E^2 = -1$ implies $ E \not\in \overline{NE}(X)_{E \geq 0}$ and $E$ is seen to generate an extremal ray in $\overline{NE}(X)$. By adjunction, $\mathcal K_X \cdot E < 0$. The contraction of the extremal ray $R = \mathbb R_{\geq 0}[E]$ is precisely the contraction of the (-1)-curve $E$. Conversely, each extremal contraction of type 1 above is the contraction of a (-1)-curve generating the extremal ray $R$. 
\end{remark}
\subsection{Equivariant contraction theorem and $G$-minimal models}
In this subsection
we prove an equivariant contraction theorem for smooth projective surfaces with finite groups of symmetries. Most steps in the proof are carried out in analogy to the proof of the standard contraction theorem.
\begin{definition}\label{equiContraction}
Let $G$ be a finite group, let $X$ be a smooth projective surface with $G$-action and let $R \subset \overline{NE}(X)^G$ be $G$-extremal ray. A morphism $\mathrm{cont}_R^G: X \to Z$ is called the \emph{$G$-equivariant contraction of $R$} if
\begin{samepage}
\begin{itemize}
\item[-] $\mathrm{cont}_R^G$ is equivariant with respect to $G$
\item[-] $(\mathrm{cont}_R^G)_*\mathcal{O}_X = \mathcal{O}_Z$ and 
\item[-] $\mathrm{cont}_R(C) = \{\text{point}\}$ for an irreducible curve $C\subset X$ if and only if $[GC] \in R$. 
\end{itemize}
\end{samepage}
\end{definition}
\begin{theorem}
Let $G$ be a finite group, let $X$ be a smooth projective surface with $G$-action and let $R$ be a $G$-extremal ray with $ \mathcal{K}_X \cdot R <0$.  Then $R$ can be spanned by the class of $C= \sum_{g \in G} gC_0$ for a rational curve $C_0$, the equivariant contraction morphism $\mathrm{cont}_R^G: X \to Z$ exists and is one of the following three types:
\begin{enumerate}
\item $C^2 <0$ and $gC_0$ are smooth disjoint (-1)-curves. The map $\mathrm{cont}_R^G: X \to Z$ is the equivariant blow down of the disjoint union $\bigcup_{g \in G} gC_0$.
\item $C^2 =0$ and any connected component of $C$ is either irreducible or the union of two (-1)-curves intersecting transversally at a single point. The map $\mathrm{cont}_R^G: X \to Z$ defines an equivariant conic bundle over a smooth curve .
\item $C^2 >0$ , $N_1(X)^G = \mathbb{R }$ and $\mathcal{K}_X^{-1}$ is ample, i.e., $X$ is a Del Pezzo surface. The map $\mathrm{cont}_R^G: X \to Z$ is constant, $Z$ is a point. 
\end{enumerate}
\end{theorem}
\begin{proof}
Let $R$ be a $G$-extremal ray with $\mathcal{K}_X \cdot R <0$.
By Lemma \ref{Gextremalray} the ray $R$ can be
 spanned by a 1-cycle of the form $C = GC_0$ for a rational curve
 $C_0$. Let $n = |GC_0|$ and write $C = \sum_{i=1}^n C_i$
 where the $C_i$ correspond to $gC_0$ for some $g \in G$.
We distinguish three cases according to the sign of the self-intersection of $C$.
\subsection*{The case $C^2 <0$}
 We write $0 > C^2 = \sum_i C_i^2 + \sum_{i\neq j}C_i \cdot C_j$. Since $C_i$ are effective curves we know $C_i \cdot C_j \geq 0$ for all $i \neq j$. Since all curves $C_i$ have the same negative self-intersection and by assumption, $\mathcal{K}_X \cdot C = \sum_i \mathcal{K}_X \cdot C_i = n(\mathcal{K}_X \cdot C_i) <0$ the adjunction formula reads
$
2g(C_i) -2 = -2 = \mathcal{K}_X  \cdot C_i + C_i^2
$. Consequently, $\mathcal{K}_X \cdot C_i = -1$ and $C_i^2 = -1$.
 It remains to show that all $C_i$ are disjoint. We assume the contrary and without loss of generality 
$C_1 \cap C_2 \neq \emptyset$. Now $gC_1 \cap g C_2 \neq \emptyset $ for all $g \in G$ and $\sum_{i \neq j} C_i \cdot C_j \geq n$. This is however contrary to $0>C^2 = \sum_i C_i^2 + \sum_{i\neq j}C_i \cdot C_j = -n +  \sum_{i\neq j}C_i \cdot C_j$.

We let $\mathrm{cont}_R^G: X \to Z$ be the blow-down of $\bigcup_{i=1}^n C_i$ which is equivariant with respect to the induced action on $Z$ and fulfills $(\mathrm{cont}_R^G)_*\mathcal{O}_X = \mathcal{O}_Z$. If $D$ is an irreducible curve such that $\mathrm{cont}_R^G(D)= \{\text{point}\}$, then $D= gC_0$ for some $g \in G$. In particular, $GD= GC_0 =\lambda C$ and $[GD] \in R$. Conversely, if $[GD] \in R$ for some irreducible curve $D$, then $[GD]= \lambda [C]$ for some $\lambda \in \mathbb R _{\geq 0}$. Now $(GD)\cdot C = \lambda C^2 <0$. It follows that $D$ is an irreducible component of $C$. 
\subsection*{The case $C^2 >0$}
This case is treated in precisely the same way as the corresponding case in the standard contraction theorem.
Our aim is to show that $[C]$ is in the interior of $\overline{NE}(X)^G$. This is a consequence of the following lemma (Corollary 1.21 in \cite{kollarmori}).
\begin{lemma}
Let $X$ be a projective surface and let $L$ be an ample line bundle on $X$. Then the set $Q = \{[E] \in N_1(X) \ |\ E^2 >0\}$ has two connected components
$Q^+=\{ [E] \in Q \ |L \cdot E >0\}$ and $Q^-=\{ [E] \in Q \ |L \cdot E <0\}$. Moreover, $Q^+ \subset \overline{NE}(X)$. 
\end{lemma}
This lemma follows from the Hodge Index Theorem and the fact that $E^2 >0$ implies that either $E$ or $-E$ is effective.

We consider an effective cycle $C = \sum C_i$ with $C^2 >0$. By the above lemma, $[C]$ is contained in $Q^+$ which is an open subset of $N_1(X)$ contained in $\overline{NE}(X)$. In particular, $[C]$ lies in the interior of $\overline{NE}(X)$. The $G$-extremal ray $R = \mathbb{R }_{\geq 0} [C]$ can only lie in the interior if $\overline{NE}(X)^G= R$. By assumption $\mathcal{K}_X \cdot R <0$, hence $\mathcal{K}_X$ is negative on $\overline{NE}(X)^G \backslash \{0\}$ and therefore on $\overline{NE}(X) \backslash \{0\}$. The anticanonical bundle $\mathcal{K}_X^{-1}$ is ample by Kleiman ampleness criterion and $X$ is a Del Pezzo surface. 

We can define a constant map $\mathrm{cont}_R^G$ mapping $X$ to a point $Z$ which is the equivariant contraction of $R = \overline{NE}(X)^G$ in the sense of Definition \ref{equiContraction}. 
\subsection*{The case $C^2 =0$}
Our aim is to show that for some $m >0$ the linear system $|mC|$ defines a conic bundle structure on $X$. The argument is seperated into a number of lemmata.
For the convenience of the reader, we include also the proofs of well-known preparatory lemmata which do not involve group actions. 
\begin{lemma}
$H^2(X, \mathcal{O}(mC)) =0$ for $m \gg 0$.
\end{lemma}
\begin{proof}
By Serre's duality, 
$
h^2(X, \mathcal{O}(mC))= h^0(\mathcal{O}(-mC)\otimes \mathcal K_X)
$.
Since $C$ is an effective divisor on $X$, it follows that $h^0(\mathcal{O}(-mC)\otimes \mathcal K_X)=0$ for $m \gg 0$.
\end{proof}
\begin{lemma}
 For $m \gg 0$ the dimension $h^0(X,\mathcal{O}(mC))$ of $H^0(X,\mathcal{O}(mC))$ is at least two.
\end{lemma}
\begin{proof}
Let $m$ be such that $h^2(X, \mathcal{O}(mC)) =0$. For a line bundle $L$ on $X$ we denote by $\chi(L) = \sum_i(-1)^i h^i(X,L)$ the Euler characteristic of $L$. Using the theorem of Riemann-Roch we obtain
\begin{align*}
h^0(X, \mathcal{O}(mC)) &\geq h^0(X, \mathcal{O}(mC))- h^1(X, \mathcal{O}(mC))\\
               & = h^0(X, \mathcal{O}(mC)) - h^1(X, \mathcal{O}(mC)) + h^2(X, \mathcal{O}(mC))\\
               & = \chi(\mathcal{O}(mC))\\
               & = \chi (\mathcal{O}) + \frac{1}{2}(\mathcal{O}(mC)\otimes\mathcal{K}_X^{-1})\cdot(mC)\\
               & \overset{C^2=0}{=} \chi (\mathcal{O}) - \frac{m}{2}\mathcal{K}_X \cdot C.
\end{align*}
Now $\mathcal{K}_XC<0$ implies the desired behaviour of $h^0(X, \mathcal{O}(mC))$. 
\end{proof}
For a divisor $D$ on $X$ we denote by $|D|$ the complete \emph{linear system of $D$}, i.e., the space of all effective divisors linearly equivalent to $D$. A point $p \in X$ is called a base point of $|D|$ if $p \in \mathrm{support}(C)$ for
all $C \in |D|$.
\begin{lemma}
 There exists $m' >0$ such that the linear system $|m'C|$ is base point free.
\end{lemma}
\begin{proof}
Let $m$ be chosen such that $h^0(X, \mathcal{O}(mC)) \geq 2$. We denote by $B$ the \emph{fixed part} of the linear system $|mC|$, i.e., 
the biggest divisor $B$ such that each $D \in |mC|$ can be decomposed as $D = B + E_D$ for some effective divisor $E_D$.
The support of $B$ is the union of all positive dimensional components of the set of base points of $|mC|$. We 
assume that $B$ is nonempty.
The choice of $m$ guarantees that $|mC|$ is not fixed, i.e., there exists $D \in |mC|$ with $D \neq B$.
Since $\mathrm{supp}(B) \subset \{s=0\}$ for all $ s \in \Gamma(X, \mathcal O(mC))$, each irreducible component of $\mathrm{supp}(B)$ is an irreducible component of $C$ and $G$-invariance of $C$ implies $G$-invariance of the fixed part of $|mC|$. It follows that $B=m_0 C$ for some $m_0 < m$.
Decomposing $|mC|$ into the fixed part $B = m_0C$  and the remaining \emph{free part} $|(m-m_0)C|$ shows that some multiple $|m'C|$ for $m' >0$  has no fixed components.
The linear system $|m'C|$ has no isolated base points since these would
correspond to isolated points of intersection of divisors linearly equivalent to $m'C$. Such intersections are excluded by $C^2=0$.
 \end{proof}
We consider the base point free linear system $|m'C|$ and the associated morphism 
$\varphi =\varphi_{|m'C|} : X \to \varphi(X) \subset \mathbb P(\Gamma(X,\mathcal{O}(m'C))^*)$.
Since $C$ is $G$-invariant, it follows that $\varphi$ is an equivariant map with respect to action of $G$ on $\mathbb P(\Gamma(X,\mathcal{O}(m'C))^*)$ induced by pullback of sections.

Let $z$ be a linear hyperplane in $\Gamma(X,\mathcal{O}(m'C))$. By definition, $\varphi^{-1}(z)= \bigcap _{s\in z}(s)_0$ where $(s)_0$ denotes the zero set of the section $s$. Since $(s)_0$ is linearly equivalent to $m'C$ and $C^2=0$, the intersection $\bigcap _{s\in z}(s)_0$ does not consist of isolated points but all $(s)_0$ with $s \in z$ have a common component. In particular, each fiber is one-dimensional.
Let $f: X \to Z$ be the Stein factorization of $\varphi: X \to \varphi(X)$. The space $Z$ is normal and 1-dimensional, i.e., $Z$ is a Riemann surface. Note that there is a $G$-action on $Z$ such that $f$ is equivariant. 
\begin{lemma}
The map $f: X \to Z$ defines an equivariant conic bundle\index{conic bundle}, i.e., an equivariant fibration with general fiber isomorphic to $\mathbb{P}_1$.
\end{lemma}
\begin{proof}
Let $F$ be a smooth fiber of $f$. By construction, $F$ is a component of $(s)_0$ for some element $s \in \Gamma(X,\mathcal{O}(m'C))$. We can find an effective 1-cycle $D$ such that $(s)_0 = F+D$. Averaging over the group $G$ we obtain
$
\sum_{g\in G}gF + \sum_{g\in G}gD = \sum_{g\in G}g(s)_0
$. 
Recalling $(s)_0 \sim m'C $ and $[C] \in \overline{NE}(X)^G$ we deduce
\[
[\sum_{g\in G}gF + \sum_{g\in G}gD] = [\sum_{g\in G}g(s)_0] = m'[\sum_{g\in G}gC]= m |G| [C].
\]
This shows that $[\sum_{g\in G}gF + \sum_{g\in G}gD]$ in contained in the $G$-extremal ray generated by $[C]$. By the definition of extremality $[\sum_{g\in G}gF] = \lambda [C] \in \mathbb{R }^{>0}[C]$ and therefore $\mathcal{K}_X \cdot (\sum_{g\in G}gF) <0$. This implies $\mathcal{K}_X \cdot F<0$.

In order to determine the self-intersection of $F$, we first observe $(\sum_{g\in G}gF)^2= \lambda^2 C^2 =0$. Since $F$ is a fiber of a $G$-equivariant fibration, we know that $\sum_{g\in G}gF = kF + kF_1 + \dots + kF_l$ where $F, F_1, \dots F_l$ are distinct fibers of $f$ and $k \in \mathbb N ^{>0}$. Now
$0=(\sum_{g\in G}gF)^2 = (l+1)k^2F^2$
shows $F^2=0$. The adjunction formula implies $g(F)=0$ and $F$ is isomorphic to $\mathbb P_1$.
\end{proof}
The map $\mathrm{cont}_R^G:=f$ is equivariant and fulfills $f_* \mathcal{O}_X = \mathcal{O}_Z$ by Stein's factorization theorem. Let $D$ be an irreducible curve in $X$ such that $f$ maps $D$ to a point, i.e., $D$ is contained in a fiber of $f$. Going through the same arguments as above one checks that $[GD] \in R$. Conversely, if $D$ is an irreducible curve in $X$ such that $[GD] \in R$ it follows that $(GD) \cdot C=0$. If $D$ is not contracted by $f$, then $f(D) = Z$ and $D$ meets every fiber of $f$. In particular, $D \cdot C >0$, a contradiction. It follows that $D$ must be contracted by $f$.

This completes the proof of the equivariant contraction theorem.
\end{proof}
A \emph{conic} is a divisor defined
by a homogeneous polynomial of degree two in $\mathbb P_2$.
It is therefore either a smooth curve of degree two and multiplicity
one, two projective lines of multiplicity one which 
intersect transversally in one point, or a single line
of multiplicity two. A smooth conic is isomorphic to $\mathbb P_1$. A conic bundle $X\to Z$ is, as
the name suggests, a ``bundle'' of conics. Its possible degenerations correspond precisely to the degenerations of conics.

The singular fibers of the conic bundle in case (2) of the theorem above are characterized by the 
following lemma stating that only conic degenerations of the first kind may occur.
\begin{lemma}\label{singular fibers of conic bundle}
Let $R = \mathbb R ^{>0} [C]$ be a $\mathcal K_X$-negative $G$-extremal ray with $C^2=0$. Let $\mathrm{cont}_R^G:=f: X \to Z$ be the equivariant contraction of $R$ defining a conic bundle structure on $X$. Then every singular fiber of $f$ is a union of two (-1)-curves intersecting transversally. 
\end{lemma}
\begin{proof}
Let $F$ be a singular fiber of $f$. The same argument as in the previous lemma yields that $\mathcal{K}_X \cdot F<0$ and $F^2 =0$. Since $F$ is connected, the arithmetic genus of $F$ is zero and $\mathcal{K}_X \cdot F = -2$. The assumption on $F$ being singular implies that $F$ must be reducible. Let $F = \sum F_i$ be the decomposition into irreducible components and note that $g(F_i)=0$ for all $i$. 

We apply the same argument as above to the component $F_i$ of $F$: after averaging over $G$ we deduce that $GF_i$ is in the $G$-extremal ray $R$ and $ \mathcal K _X \cdot F_i <0$.
Since $-2 = \mathcal K _X \cdot F = \sum \mathcal K _X \cdot F_i$, we may conclude that $F = F_1 +F_2$ and $F_i^2=-1$. 
The desired result follows.
\end{proof}
It should be remarked that
 $G$-equivariant conic bundles with or without singular fibers can be studied by considering the $G$-action on the base and the actions of the isotropy groups of points of the base on the corresponding fibers. 
\subsection {$G$-minimal models of surfaces}
Let $X$ be a surface with an action of a finite group $G$ such that $\mathcal{K}_X$ is not nef, i.e., $\overline{NE}(X)_{\mathcal{K}_X<0}$ is nonempty. 
\begin{lemma}
There exists a $G$-extremal ray $R$ such that $\mathcal{K}_X \cdot R<0$. 
\end{lemma}
\begin{proof}
Let $[C] \in \overline{NE}(X)_{\mathcal{K}_X<0}\neq \emptyset$ and consider $[GC] \in \overline{NE}(X)^G$. The $G$-orbit or $G$-average of a $\mathcal{K}_X$-negative effective curve is again $\mathcal{K}_X$-negative. It follows that $\overline{NE}(X)^G_{\mathcal{K}_X<0}$ is nonempty. 
Let $L$ be a $G$-invariant ample line bundle on $X$.
By the cone theorem,
for any $\varepsilon>0$
\begin{equation}\label{conethmformula}
\overline{NE}(X)^G = \overline{NE}(X)^G _{(\mathcal{K}_X+\varepsilon L) \geq 0} + \sum_{\text{finite}} \mathbb R_{\geq 0} G[C_i]. 
\end{equation}
with $\mathcal K _X \cdot C_i < 0$ for all $i$.
Since $\overline{NE}(X)^G_{\mathcal{K}_X<0}$ is nonempty, we may choose $\varepsilon>0$ such that $\overline{NE}(X)^G \neq \overline{NE}(X)^G _{(\mathcal{K}_X+\varepsilon L) \geq 0}$.
If the ray $R_1 = \mathbb R_{\geq 0} G[C_1]$ is not extremal in $\overline{NE}(X)^G$, then its generator $G[C_1]$ can be decomposed as a sum of elements of $\overline{NE}(X)^G$ not contained in $R_1$. It follows that 
the ray $R_1$ is superfluous in the formula \eqref{conethmformula}.
Since
$\overline{NE}(X)^G \neq \overline{NE}(X)^G _{(\mathcal{K}_X+\varepsilon L) \geq 0}$ by assumption, we may not remove all rays $R_i$ from the formula and at least one ray $R_i = \mathbb R_{\geq 0} G[C_i]$ is $G$-extremal.
\end{proof}
We apply the equivariant contraction theorem to $X$:
In case (1) we obtain a new surface $Z$ from $X$ by blowing down a $G$-orbit of disjoint (-1)-curves. There is a canonically defined holomorphic $G$-action on $Z$ such that the blow-down is equivariant. If $K_Z$ is not nef, we repeat the procedure which will stop after a finite number of steps. In case (2) we obtain an equivariant conic bundle structure on $X$. In case (3) we conclude that $X$ is a Del Pezzo surface with $G$-action. We call the $G$-surface obtained from $X$ at the end of this procedure a \emph{$G$-minimal model of $X$}.

As a special case, we consider a rational surface $X$ with $G$-action. 
Since the canonical bundle $\mathcal{K}_X$ of a rational surface $X$ is never nef, 
a $G$-minimal model of $X$ is an equivariant conic bundle over a smooth rational curve $Z$ or a Del Pezzo surface with $G$-action. 
This proves the well-known classification of $G$-minimal models of rational surfaces (cf.\,\cite{maninminimal}, \cite{isk}).

Although this classification does classically not rely on Mori theory, the proof given above is based on Mori's approach. We therefore refer to an equivariant reduction $Y \to Y_\mathrm{min}$ as an \emph{equivariant Mori reduction}.
%
%
%
%
\section {Del Pezzo surfaces}
The equivariant minimal model program for surfaces presents us with the task of studying automorphism groups of Del Pezzo surfaces. 

A Del Pezzo surfaces is defined as a connected compact complex surface
whose anticanonical line bundle is ample. Using Kleiman's ampleness
criterion we have identified the class of Del Pezzo surfaces as a
class of $G$-minimal surfaces. Here we wish to replace this very
abstract notion of ampleness by a definition involving the notions of bundles and sections. 

We let $X$ be a connected compact complex surface, $T_X$ the holomorphic tangent bundle,
$T_X^*$ the contangent bundle, and $\mathcal K_X:=\Lambda ^2T^*_X$
its top exterior power, the canonical line bundle. The 
\emph{anticanonical bundle} is given by $\mathcal K^{-1}_X=\Lambda ^2T_X$.

In this terminology the space
$V_k$, which we have discussed in a naive fashion up to this
point, is the space of sections $\Gamma (X,\mathcal K_X ^k)$. 
The requirement that $\mathcal K^{-1}_X$ is \emph{ample} means that some power $(\mathcal K_X^{-1})^k =:-k \mathcal K_X$
has many sections.  More precisely, one requires that
for some $k$ the map $X\to \mathbb P(\Gamma (X,-k \mathcal K_X)^*)$
is a holomorphic embedding. 

Using wedge-products of holomorphic vector fields one
easily shows that the anticanonical bundles of $\mathbb P_2$
and $\mathbb P_1\times \mathbb P_1$ are ample. In order to characterize
the remaining
Del Pezzo surfaces it is convenient to introduce the \emph{degree}
of a Del Pezzo surface as the self-intersection number $d$ 
of an anticanonical divisor. 
For $\mathbb P_2$ it is simple to
compute this degree:  $\mathcal K^{-1}_{\mathbb P_2}$ is the 3rd power $H^3$
of the hyperplane bundle and its sections are homogeneous
polynomials of degree three. By Bezout's theorem, the intersection
of two such cubics consists of nine points counted with multiplicity
and the degree $d$ of the Del Pezzo surface $\mathbb P_2$ equals $d=9$. The possible degrees of Del Pezzo surfaces range from one to nine.

The following theorem (cf.\,Theorem 24.4 in \cite{M}) gives a classification
of Del Pezzo surfaces according to their degree.
\begin{theorem}
Let $Z$ be a Del Pezzo surface of degree $d$. 
\begin{samepage}
\begin{itemize}
\item[-]
If $d=9$, then $Z$ is isomorphic to $\mathbb P_2$.
\item[-]
If $d=8$, then $Z$ is isomorphic to either $\mathbb P_1 \times \mathbb
P_1$ or the blow-up of $\mathbb P_2$ in one point.
\item[-]
If $1 \leq d \leq 7$, then $Z$ is isomorphic to the blow-up of
$\mathbb P_2$ in $9-d$ points in general position, i.e., no three
points lie on one line and no six points lie on one conic. 
\end{itemize}
\end{samepage}
\end{theorem}
Using the theorem above we may identify a Del Pezzo surface $X \neq \mathbb P_1 \times \mathbb P_1$ with the blow up $X_{\{p_1,\ldots ,p_m\}}$ of $\mathbb P_2$ at each point of $\{p_1,\ldots ,p_m\}$ 
for $ m \in \{0, \dots, 8\}$.
We carry the points $p_1,\ldots ,p_m$  in the notation, because for $m\ge 5$
the complex structure really does depend on the points which 
are blown up.  For example, if $m=4$, using automorphisms
of $\mathbb P_2$ we can move the points to a standard 
location, e.g.,\,$[1:0:0],[0:1:0],[0:0:1],[1:1:1]$.  Of
course such a normal form is even easier to achieve if
$m<4$.  On the other hand, if $m>4$, then we put 
$p_1,\ldots, p_4$ in this normal form, but
$p_5,\ldots, p_m$ are allowed to vary.  As these points
vary the complex structure of the Del Pezzo surface
varies.  So for $m\ge 5$ Del Pezzo surfaces come in families.
\subsection{Automorphism groups of Del Pezzo surfaces}
Let us now turn to a study of automorphism groups of Del Pezzo 
surfaces.  
In analogy to the case of surfaces of general type we 
begin with the following first remark.
\begin {proposition}
The automorphism group $\mathrm {Aut}(X)$ of a Del Pezzo surface
is an algebraic group acting algebraically on $X$.
\end {proposition}
\begin {proof}
For $k$ sufficiently large the map 
$X\to \mathbb P(\Gamma (X,-k \mathcal K_X)^*)$ is an 
$\mathrm {Aut}(X)$-equi\-va\-riant, holomorphic embedding. We denote its image by $Z$.
The group $\mathrm {Aut}(X)$ can be realized as the stabilizer
of $Z$ in the automorphism group of the ambient 
projective space, which is itself an algebraic group. 
By the theorem of Chow the space $Z$ is an algebraic subvariety.
Since the stabilizer of an algebraic subvariety is
an algebraic group acting algebraically,
the result follows.
\end {proof}
There are Del Pezzo surfaces with positive-dimensional
automorphism groups, e.g., $\mathbb P_2$ and 
$\mathbb P_1\times \mathbb P_1$ are even homogeneous.
Since the stabilizer in $\mathrm{Aut}( \mathbb P_2)$ of $\{p_1,\ldots ,p_m\}$ has an open orbit for $m \leq 3$, the following shows in particular
that the Del Pezzo surfaces of degree at least six are
almost homogeneous, i.e., their  
automorphism groups have an open orbit.
\begin {proposition}\label{lifting}
If $\pi :X_{\{p_1,\ldots ,p_m\}}\to \mathbb P_2$ is the defining
blow up of the Del Pezzo surface $X_{\{p_1,\ldots ,p_m\}}$ and
$g\in \mathrm {Aut}(\mathbb P_2)$ stabilizes the
set $\{p_1,\ldots ,p_m\}$, then there exists a uniquely
defined automorphism $\hat g\in \mathrm {Aut}(X_{\{p_1,\ldots ,p_m\}})$ so that
$g\circ \pi =\pi \circ \hat g$.
\end {proposition}
\begin {proof}
If $E_j$ is the $\pi $-preimage of $p_j$, $j=1,\ldots ,m$,
then there exists a unique automorphism $\hat g$ of
$X_{\{p_1,\ldots ,p_m\}}\setminus \bigcup_j E_j$ with the desired
property. Thus it is a matter of extending $\hat g$
to the full Del Pezzo surface. It is enough to show
that it extends across $E=E_1$. 
For notational simplicity we may assume
that $gp_1=p_1=:p$. Since points of $E$
correspond to tangent lines through $\pi(E) = p$ and $gp=p$,
it follows that $\hat g$ extends continuously across $E$.
The fact that this extension is holomorphic is guaranteed
by Riemann's Hebbarkeitssatz. It is
an automorphism since $g^{-1}$ can likewise be lifted to $X_{\{p_1,\ldots ,p_m\}}$.
\end {proof}
Furthermore, the considerations in this section will benefit from the following converse of the above statement.
\begin {proposition}\label {small automorphisms}
The defining projection $\pi :X_{\{p_1,\ldots ,p_m\}}\to \mathbb P_2$
is equivariant with respect to the connected component
$\mathrm {Aut}(X_{\{p_1,\ldots ,p_m\}} )^\circ $ containing the identity.
\end {proposition}
\begin {proof}
Let $g_t$ be a a 1-parameter subgroup in 
$\mathrm {Aut}(X_{\{p_1,\ldots ,p_m\}})$ and let $E:=E_j$ be the $\pi $-preimage of
$p:=p_j$. Given a small neighborhood $U$ of $p$, if
$t$ is sufficiently small, then $\pi(g_t E)$ is contained in $U$. 
Since $\pi $ is a proper holomorphic map, it follows from Remmert's mapping theorem that
$\pi(g_t E)$ is a compact analytic subset of $U$. If $U$
is sufficiently small, its only analytic subsets are discrete. Hence, $\pi(g_t E)$ is a single point and $g_t E=E$.  
As a result  
every 1-parameter subgroup of $\mathrm {Aut}(X_{\{p_1,\ldots ,p_m\}})$ 
stabilizes every $E_j$ and consequently the same is true for the
full connected component $\mathrm {Aut}(X_{\{p_1,\ldots ,p_m\}})^\circ $.   
Thus every element of $\mathrm {Aut}(X_{\{p_1,\ldots ,p_m\}})^\circ $
defines a map of $\mathbb P_2$ fixing $\{p_1,\ldots ,p_m\}$
which is holomorphic outside $\{p_1,\ldots ,p_m\}$ and the desired result follows
from Riemann's Hebbarkeitssatz.
\end {proof}
\begin {corollary}
If $m\ge 4$, then the automorphism group of the Del Pezzo
surface $X_{\{p_1,\ldots ,p_m\}}$ is finite.
\end {corollary}
Detailed lists of the automorphism groups Del Pezzo surfaces
can be found in Dolgachev's book (\cite {D}).  For the convenience
of the reader we present a brief road map here.
\subsection* {The projective plane}
The action of $\mathrm {SL}_3(\mathbb C)$ on $\mathbb C^3$
defines a surjective homomorphism 
$\mathrm {SL}_3(\mathbb C) \to \mathrm {Aut}(\mathbb P_2)$ and $\mathrm {Aut}(\mathbb P_2)$ can be identified with
the quotient $\mathrm {SL}_3(\mathbb C)/C_3$. Here $C_3$
is embedded as the center of $\mathrm {SL}_3(\mathbb C)$ which
consists of matrices of the form $\lambda \mathrm {Id}$ with
$\lambda ^3=1$.  
In general, the action of a finite group $G$ on $\mathbb P_2$ is given by a 3-dimensional linear representation of a nontrivial central extension of $G$. An interesting example is provided by the Valentiner group defining an action of the alternating group $A_6$ on $\mathbb P_2$.
\subsection* {The space $\mathbb P_1\times \mathbb P_1$}
The  group
$S:=\mathrm {Aut}(\mathbb P_1)\times \mathrm {Aut}(\mathbb P_1)$
acts on the Del Pezzo surface $ X = {\mathbb P_1\times \mathbb P_1}$ in an obvious fashion.
  In addition, $\mathrm{Aut}(X)$ contains the holomorphic involution $\sigma $ which exchanges
the factors.
The group $S$ coincides with the subgroup of all $h \in \mathrm{Aut}(X)$ such that both projections are $h$-equivariant. At the level of homology classes $F_i$ of fibers, $hF_i = F_i$ for each $h \in S$. 
In fact $S$ can be identified with the subgroup of elements of $\mathrm {Aut}(X)$ having this
property. For an arbitrary $ g\in \mathrm{Aut}(X)$, since $(gF_i)^2=0$, either $gF_i=F_i$ for $i=1,2$ or $gF_1=F_2$. So either $g \in S$ or $\sigma \circ g \in S$  
and 
$\mathrm {Aut}(\mathbb P_1\times \mathbb P_1)=
S\rtimes \langle \sigma \rangle $.

\medskip
In the following we denote by $X_d$ a Del Pezzo
surface of degree $d$ which arises as the blow up of
$\mathbb P_2$ in $9-d$ points. We have already pointed out that any Del Pezzo 
surface except ${\mathbb P_1\times \mathbb P_1}$ is of this form.
\subsection* {The simple blow up of $\mathbb P_2$}
The surface $X_8$
is the simple blow up $\mathrm {Bl}_p(\mathbb P_2)$.
Since the exceptional curve $E$ of this blow up
is the unique (-1)-curve in $X_8$,
it follows that every automorphism of $X_8$ stabilizes
$E$ and $\mathrm {Aut}(X_8)$ can be identified with
the isotropy group at $p$ in $\mathrm {Aut}(\mathbb P_2)$
\subsection* {The blow up of ${\mathbb P_2}$ in two points}
Let $p_1,p_2\in \mathbb P_2$ be the two points which
are blown up to obtain $X_7$. Note that there exists
an automorphism $\sigma $ of $\mathbb P_2$, more precisely a holomorphic
involution,  which exchanges
$p_1$ and $p_2$, and note that by the arguments given above
it can be lifted to an automorphism of $X_2$. Denoting
this lift by $\hat \sigma $, it follows that
$\mathrm {Aut}(X_2)$ can be identified with the semidirect product 
$T\rtimes \langle \hat\sigma \rangle$ where $T$ is the subgroup
of $\mathrm {Aut}(\mathbb P_2)$ consisting of transformations
fixing both $p_1$ and $p_2$.
\subsection* {The three-point blow up of ${\mathbb P_2}$}
Let $X_6$ be the Del Pezzo surface defined by blowing up
the three points $p_1=[1:0:0]$, $p_2=[0:1:0]$, and $p_3=[0:0:1]$.  
As we know, the connected component of the identity in 
$\mathrm {Aut}(X_6)$ is the connected component
of the stabilizer of $\{p_1,p_2,p_3\}$.  This is the group
of diagonal matrices in $\mathrm {SL}_3(\mathbb C)$
modulo the center of $\mathrm {SL}_3(\mathbb C)$ and is therefore
isomorphic to $(\mathbb C^*)^2$.  The full permutation group
of $\{p_1,p_2,p_3\}$ can also be realized in 
$\mathrm {Aut}(\mathbb P_2)$. We see that
the subgroup of automorphisms of $X_6$ which are equivariant
with respect to the defining map $X_6\to \mathbb P_2$ is  
isomorphic to $(\mathbb C^*)^2\rtimes S_3$.

The proper transform $\hat L_{ij}$ in $X_6$ of the projective line $L_{ij}$ joining $p_i$ and
$p_j$ is a (-1)-curve. This is due to the fact that $L_{ij}^2=1$ 
and $L_{ij}$ contains exactly two points which are blown up. 
One can show that the only (-1)-curves in $X_6$ are the exceptional curves $E_i$ obtained from blowing up the points $p_i$ and the ``lines''  $\hat L_{ij}$. The graph of this configuration of curves is a hexagon $H$. Restriction yields a 
homomorphism $R: \mathrm {Aut}(X_6)\to I(H)$, where
$I(H) \cong D_{12}$ is the group of rigid motions of the hexagon and the kernel of 
$R$ is the connected component $\mathrm {Aut}(X_6)^\circ \cong (\mathbb C^*)^2$ discussed above. 

We have already seen that the permutation group $S_3$ is
contained in the image of $R$ and will now show that
there is an additional involution in this image so that
in fact $R$ is surjective.  In order to determine this
involution it is useful to regard $X_6$ as the blow up
of $\mathbb P_1\times \mathbb P_1$ in the 
\emph{antidiagonal corner points} $c_1=([1:0],[0:1])$ and
$c_2=([0:1],[1:0])$.  If one draws 
$\mathbb P_1\times \mathbb P_1$ as a square, then the
hexagon $H$ consists of the (-1)-curves arising from 
blowing up $c_1$ and $c_2$ together with the proper transforms
of the four edges of the square.  Using elementary
intersection arguments one can show that by blowing
down either of the two configurations of three disjoint 
curves in this hexagon one obtains $\mathbb P_2$ with three points in general position. Hence,
this blow up is indeed $X_6$ and the additional involution discussed
above is defined by lifting the involution of 
$\mathbb P_1\times \mathbb P_1$ which exchanges the factors.
Thus we have shown that $\mathrm {Aut}(X_6)$ is naturally
isomorphic to $(\mathbb C^*)^2\rtimes I(H)=(\mathbb C^*)^2\rtimes D_{12}$.
\subsection* {The Del Pezzo surface of degree five}
We may define $X_5$ by blowing up $p_4:=[1:1:1]$ and $p_1$, $p_2$, $p_3$ as
above. Defining $\hat L_{ij}$
as before we obtain a configuration of ten (-1)-curves which 
in fact is the entire collection of (-1)-curves in $X_5$.  The
dual graph of this configuration is known as the Petersen graph $P$.
\begin{figure}[h]
\begin{center}
\includegraphics[width = 0.2\textwidth, angle=180]{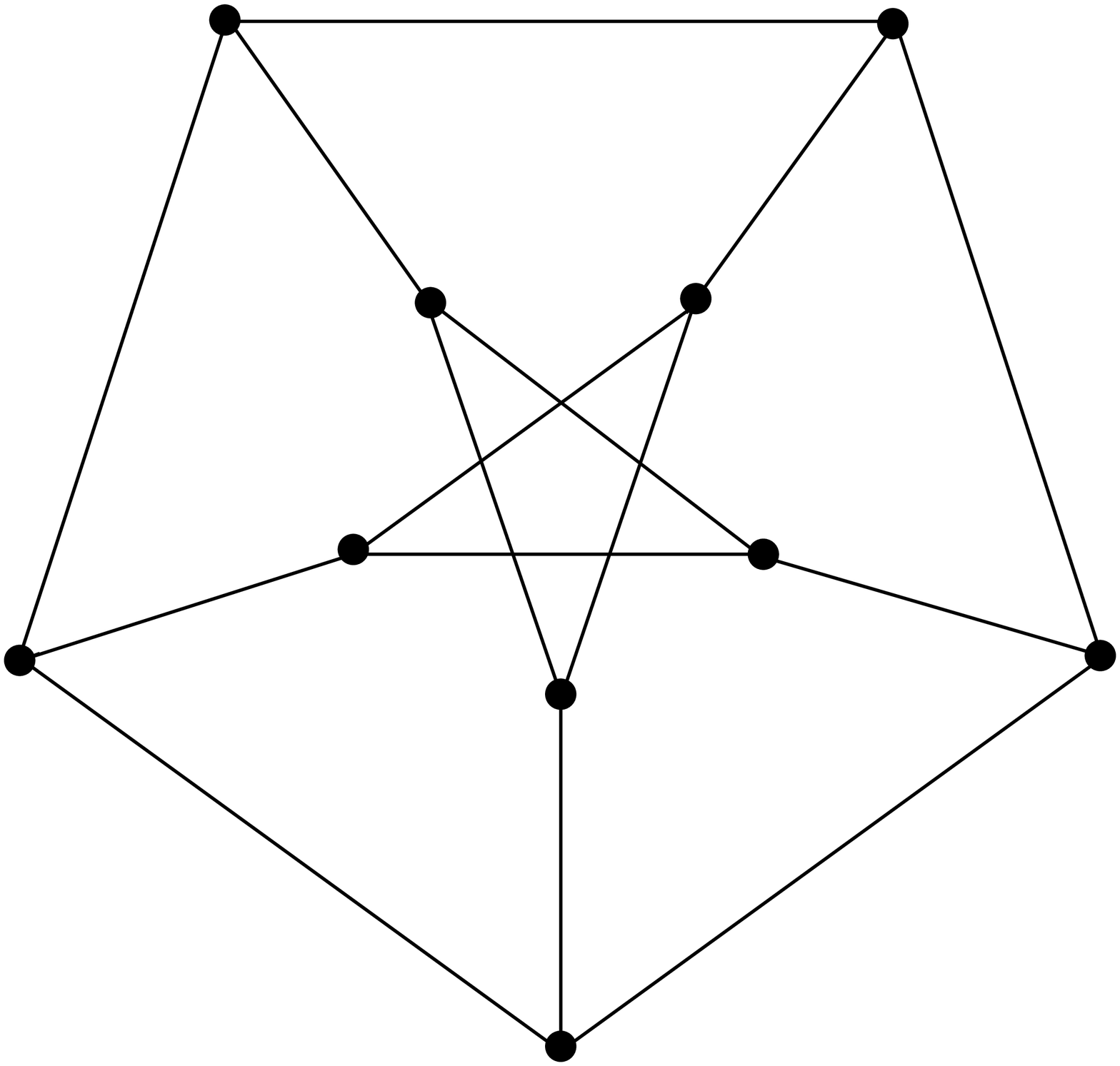}
\caption{The Petersen graph}
\end{center}
\end{figure}

A point in the graph represents a (-1)-curve and
two such curves intersect (trans\-ver\-sal\-ly) if and only 
if the corresponding points in the graph are connected by a line segment. The connected 
component of $\mathrm {Aut}(X_5)$ is trivial and therefore
the restriction map $R:\mathrm {Aut}(X_5)\to I(P)$ realizes
the automorphism group of $X_5$ as a subgroup of the graph
automorphism group $I(P)\cong S_5$.
 
There are various collections of four disjoint curves in $P$
which can be blown down to obtain a copy of $\mathbb P_2$ with four distinguished points. By Proposition \ref{lifting} their permutation group $S_4$ can be identified with a subgroup of $\mathrm{Aut}(X_5)$.
For two configurations of disjoint curves in $P$ which have one curve in common
we observe that the two corresponding copies of $S_4$
do not define the same subgroup of $\mathrm {Aut}(X_5)$.
Thus together they generate a subgroup of $S_5$
properly containing $S_4$. Since the index of $S_4$
in $S_5$ is prime, it follows that they generate the
full group $S_5$. We have shown  $\mathrm {Aut}(X_5)\cong I(P)\cong S_5$.

\medskip
In the remaining cases we use a number of basic general
facts about Del Pezzo surfaces. For their proofs, we refer the reader, e.g., to \cite {M} and \cite {D}.
Here we present outlines of the arguments required to identify $\mathrm{Aut}(X_d)$ for $ d \leq 4$.
\subsection* {Five-point blow ups of ${\mathbb P_2}$}
We define a surface $X_4$ by blowing up the four points $p_1, \dots, p_4$ as
above and in addition a fifth point $p_5$.  As $p_5$ moves
so does the complex structure of $X_4$.  Thus, it is to be
expected that the automorphisms group of $X_4$ depend on the position of $p_5$.

In the following, we strongly use the fact that $X_4$ is embedded by
$\mathcal K^{-1}_X$ in $\mathbb P_4$ as a surface which is the
transversal intersection of two nondegenerate quadric 3-folds.
Choosing coordinates appropriately we may assume that
these quadrics are defined by $Q_1:=\sum z_j^2$ and 
$Q_2:=\sum a_jz_j^2$ with
$a_i\not=a_j$ for $i\not=j$. Since the embedding in $\mathbb P_4$
is $\mathrm {Aut}(X_4)$-equivariant, the group 
$\mathrm {Aut}(X_4)$ can be identified with the stabilizer $S$
in $\mathrm {Aut}(\mathbb P_4)$ of the subspace 
$V:=\mathrm {Span}(Q_1,Q_2)$ in the space of all quadratic
forms. Computing in $\mathrm {SL}_5(\mathbb C)$ one sees
that $S$ is the normalizer of the group of diagonal matrices $T$
modulo the center $C_5$ of $\mathrm {SL}_5(\mathbb C)$.
This is the group $T\rtimes S_5$ where $S_5$ is acting by 
permuting the coordinate functions $z_0,\ldots ,z_4$.

The meromorphic map $\mathbb P_4\to \mathbb P(V)$ defined by $Q_1,Q_2$ is 
$S$-equivariant and therefore defines a homomorphism
$S\to S/I\hookrightarrow \mathrm {Aut}(\mathbb P_1)$. The kernel $I$
consists of those transformations which act on $Q_1$ and
$Q_2$ by the same character. Since $a_i\not=a_j$ for 
$i\not =j$, it follows that $I\cong C_2^4$ is generated
by the elements of $T$ of the form 
$\mathrm {Diag}(\pm 1,\ldots ,\pm 1)$.
Since $S_5$ normalizes this group, we see that 
$S=C_2^4\rtimes S/I$, where $S/I$ is on the one hand a subgroup
of $S_5$ and on the other a subgroup of $\mathrm {Aut}(\mathbb P_1)$.
Using this information one can directly compute all possibilities 
for $\mathrm {Aut}(X_4)$, namely $S/I \in \{C_2,C_4, S_3, D_{10} \}$. 
(see $\S$ 10.2.2 in \cite {D}).
\subsection* {Cubic surfaces}
The case of Del Pezzo surfaces of degree three is conceptually
simple, but computationally complicated.  In this
case $\mathcal K^{-1}_X$ is still very ample and embeds $X_3$ as a cubic
surface in $\mathbb P_3$, i.e., as the zero-set of a cubic
polynomial $P_3$. Since this embedding is 
$\mathrm {Aut}(X_3)$-equivariant, $\mathrm {Aut}(X_3)$
can be identified with the stabilizer in 
$\mathrm {Aut}(\mathbb P_3)$ of the line $\mathbb C \cdot P_3$
in the space of all cubics. Consequently, the classification of 
automorphism groups of Del Pezzo surfaces of degree three
amounts to the determination of the invariants of 
actions of the finite subgroups of $\mathrm {SL}_4(\mathbb C)$ 
on the space of cubic homogeneous polynomials. This is
carried out in \cite {D} where the results are presented in Table 10.3.  
\subsection* {Double covers ramified over a quartic}
A Del Pezzo surface $X_2$ of degree two,
can be realized as a 2:1 cover ramified over a smooth curve $C$ of degree
four by the anticanonical map 
$\varphi_{\mathcal K^{-1}_X}:X_2\to \mathbb P_2$.
Conversely, if $X\to \mathbb P_2$ is a 2:1 cover
ramified over a smooth quartic curve, then $X$ is
a Del Pezzo surface of degree two.

A smooth quartic curve $C$ is abstractly a Riemann surface
of genus three which is embedded as a quartic curve in $\mathbb P_2$
by its canonical bundle. This embedding is
equivariant and consequently $\mathrm {Aut}(C)$ is
the stabilizer of $C$ in $\mathrm {Aut}(\mathbb P_2)$.
Furthermore, 
$\mathrm {Aut}(C)$ is acting canonically on the bundle
space $H$ of the hyperplane section bundle because its restriction to $C$ is the canonical line bundle. 

The 2:1 cover $X_2$ of $\mathbb P_2$ ramified along $C$
is constructed in the bundle space $H^2$ where 
$\mathrm {Aut}(C)$ also acts.  So on the one hand,
the equivariant covering map $X_2\to \mathbb P_2$
defines a homomorphism of $\mathrm {Aut}(X_2)$ onto a 
subgroup of $\mathrm {Aut}(C)$, and on the other hand,
$\mathrm {Aut}(C)$ lifts to a subgroup of $\mathrm {Aut}(X_2)$.
Since the kernel of the surjective homomorphism $\mathrm {Aut}(X_2)\to \mathrm {Aut}(C)$
is generated by the covering transformation, it follows
that we have a canonical splitting 
$\mathrm {Aut}(X_2)=\mathrm {Aut}(C)\times C_2$.  
Since $C$ is equivariantly embedded in $\mathbb P_2$ and
$\mathrm {Aut}(C)$ is acting as a subgroup of 
$\mathrm {SL}_3(\mathbb C)$, the classification of the
automorphism groups of Del Pezzo surfaces of degree two
results from the classification of the finite subgroups
of $\mathrm {SL}_3(\mathbb C)$ (\cite {B,Y})
and the invariant theory of their representations on 
the space of homogeneous polynomials of degree four
(see Table 10.4 in \cite {D}).
\subsection* {Del Pezzo surfaces of degree one}
In this case the anticanonical
map is a meromorphic map $\varphi _{\mathcal K^{-1}_X}:X_1\to \mathbb P_1$
with exactly one point $p$ of indeterminacy, a so-called base point.  Thus
$p$ is fixed by $\mathrm {Aut}(X_1)$ and, since it is
a finite group and in particular compact, the linearization
of $\mathrm {Aut}(X_1)$ on $T_pX_1$ is a faithful representation.
This already places a strong limitation on the group $\mathrm {Aut}(X_1)$.
Furthermore, the map $X_1\to \mathbb P_3$ defined by
$-2 \mathcal K_X$ has no points of indeterminacy and realizes
$X_1$ as a 2:1 ramified cover over a quadric cone. A study of
these two maps, both of which are equivariant, leads to
a precise description of all possible automorphism groups
of Del Pezzo surface of degree one (Table 10.5 in \cite {D})
%
%
\section {K3-surfaces with special symmety}\label {involutions}
In this section a setting is considered where the equivariant minimal model 
program has been implemented to prove classification theorems
for K3-surfaces with finite symmetry groups.
Additional techniques which aid in simplifying the combinatorial
geometry involved in the Mori reduction are outlined and recent
results which appear in \cite {Fra} are sketched. Details of
a concrete situation involving the group $A_6$ are given in
next section.  
\subsection {Maximal groups}
A \emph{K3-surface} $X$ is a simply-connected compact complex 
surface admitting a globally defined nowhere-vanishing holomorphic 2-form $\omega$.
A transformation $g\in \mathrm {Aut}(X)$ is said to be
\emph{symplectic} if $g^*\omega =\omega $ and the
group of symplectic automorphisms is denoted by $\mathrm {Aut}_\mathrm{sym}(X)$.
If $\chi :\mathrm {Aut}(X)\to\mathbb C^*$ denotes the character
defined by $g^*\omega =\chi (g)\omega $, then
$\mathrm {Aut}_{\mathrm {sym}}(X)=\mathrm{Ker}(\chi )$.  For a finite
subgroup $G\subset \mathrm {Aut}(X)$ we have the exact
sequence 
$$
1\to G_{\mathrm {sym}}\to G\to C_n\to 1\,,
$$
where the homomorphism $G\to C_n$ is the restriction 
$\chi \vert_G$.  
The group $G$ can be regarded as a coextension of $G_{\mathrm {sym}}$
by $C_n$. Although we restrict here to the case where
$G$ is finite, it should be underlined that the full
group $\mathrm {Aut}(X)$ may not be finite.
\begin{example}
Let $T_1$ be the 1-dimensional torus defined by the lattice 
$\langle 1,i\rangle _{\mathbb Z}$ and 
let $T=T_1\times T_1=\mathbb C^2/\Lambda $. The group
$\Gamma:=\mathrm {SL}_2(\mathrm Z)$ is contained in
$\mathrm {Aut}(T)$ and centralizes the involution
$\sigma :=-\mathrm {Id}$.  Thus $\Gamma $ acts
as a group of holomorphic transformations on the quotient
$Y:=T/\sigma $. The set of singular points in $Y$
consists of 16 ordinary 
double points. The desingularization 
$\mathrm {Kum}(T)\to Y$ blows up each of these points, replacing
them by copies $E$ of $\mathbb P_1$ with $E \cdot E=-2$. The
group $\Gamma $ lifts to act as a group of holomorphic
transformations on the \emph{Kummer surface} $\mathrm {Kum}(T)$.
The holomorphic 2-form $dz\wedge dw$ on $T$ is $\sigma $-invariant
and defines a nowhere vanishing holomorphic 2-form on 
$\mathrm {Kum}(T)$ which is $\Gamma $-invariant. 
Since $X=\mathrm {Kum}(T)$ is simply-connected,
it follows that it is a K3-surface with 
$\Gamma \subset \mathrm {Aut}_{\mathrm {sym}}(X)$.
\end{example}
We are interested in finite subgroups $G$ of $\mathrm {Aut}(X)$
where $G_{\mathrm {sym}}$ is either large or possesses interesting
group structure and the following theorem of
Mukai (\cite {Mu}) is of particular relevance (See \cite {Kondo}
for an alternative proof.).
\begin {theorem} 
If $G_{\mathrm {sym}}$ is a finite group of
symplectic transformations of a K3-surface $X$, then it is contained in one
of the groups $M$ listed in the following table.
\end {theorem}
\renewcommand{\baselinestretch}{1.25}
\begin{table}[h]
\centering
\begin{tabular}{l|l|l|l}
  & $M$ & $|M|$ & \textbf{Structure} \\ \hline 
1 & $L_2(7)$ & 168 & $\mathrm {PSL}_2(\mathbb F_7)=\mathrm {GL}_3(\mathbb F_2)$\\ \hline
2 & $A_6$ & 360 & $\mathrm {even \ permutations}$\\ \hline  
3 & $S_5$ & 120 & $A_6\rtimes C_2$\\ \hline 
4 & $M_{20}$ & 960 & $C_2^4\rtimes A_5$\\ \hline
5 & $F_{384}$ & 384 & $C_2^4\rtimes S_4$\\ \hline
6 & $A_{4,4}$ & 288 & $C_2^4\rtimes A_{3,3}$\\ \hline 
7 & $T_{192}$ & 192 & $(Q_8*Q_8)\rtimes S_3$\\ \hline 
8 & $H_{192}$ & 192 & $C_2^4\rtimes D_{12}$\\ \hline 
9 & $N_{72}$ & 72 & $C_3^2\rtimes D_8$\\ \hline 
10 & $M_9$ & 72 & $C_3^2\rtimes Q_8$\\ \hline 
11 & $ T_{48}$ & 48 & $Q_8\rtimes S_3$
\end{tabular}
\end{table}
\renewcommand{\baselinestretch}{1.1}
For brevity we refer to the groups which are listed in this
table as \emph{Mukai groups}.
It should be mentioned that this result is sharp in two
senses.  First, given two \emph{maximal} groups $M_1$ and $M_2$ 
which are listed in this table, $M_1$ can not be realized as a subgroup
of $M_2$ and vice versa. Secondly, given a group $M$ listed
in the table, there exists a K3-surface $X$ with
$M\subset \mathrm {Aut}_{\mathrm {sym}}(X)$. 
\subsection {K3-surfaces with antisympletic involutions}
An element $\sigma \in \mathrm {Aut}(X)$ of order two
with $\sigma ^*\omega =-\omega $ is called an
\emph {antisymplectic involution}.  For various reasons
K3-surfaces equip\-ped with such involutions are of
particular interest, e.g., from the point of view
of moduli spaces and associated automorphic forms
(see \cite {Yo}). 

Our study of these surfaces was
motivated by an attempt to understand the K3-surfaces
which possess finite groups $G$ of automorphism where
$G_{\mathrm {sym}}$ is large, e.g., where $G_{\mathrm {sym}}$
is maximal in the sense of Mukai's Theorem.  In that
setting there are strong restrictions on the
group structure of the coextension
$1\to G_{\mathrm {sym}}\to G\to C_n\to 1$ and the
size of $n$ which show that
understanding the case $n=2$ is of particular importance.
Thus, as a starting point, we undertook the classification
project in the case where $G=G_{\mathrm {sym}}\times C_2$
and where the antisympletic involution $\sigma $ which
generates $C_2$ has a nonempty set of fixed points. 

Before turning to an outline
of the main results of \cite {Fra} we would like to
emphasize that our work was motivated by a number
very interesting works of Keum, Oguiso, Zhang
(see \cite {OZ},\cite {KOZLeech,KOZExten}) and of 
course depends on the foundational results
of Nikulin (\cite {NikulinFinite} and Mukai (\cite {Mu}).

Simplifying the notation, we are interested in classifying
triples $(X,G,\sigma)$ where $X$ is a K3-surface on
which $G$ is acting as a group of symplectic automorphisms
centralized by an antisymplectic involution $\sigma $.
We assume that $G$ is acting effectively, i.e., that the
only element of $G$ which fixes every point of $X$ is the
identity, and we wish to classify these triples up to
equivariant isomorphism.  The fixed point set
$\mathrm {Fix}(\sigma )$ is either empty or 1-dimensional
and $G$ acts naturally on the quotient $Y:=Y/\sigma $.
If $\mathrm {Fix}(X)=\emptyset $, then $Y$ is an Enriques
surface.   

If $\mathrm {Fix}(\sigma )\not=\emptyset $, the quotient
$Y$ is a smooth rational surface. 
Thinking in terms of the method of 
quotients by small subgroups which was discussed in $\S$ \ref{smallgroups},
we have moved to a $G$-manifold $Y$ of lower Kodaira-dimension. In this case we
apply the equivariant minimal model program to obtain an equivariant Mori reduction $Y\to Y_1\to \ldots \to Y_N=Y_\mathrm{min}$.
It follows that  $Y_\mathrm{min}$ is either a Del Pezzo surface or an 
equivariant conic bundle over $\mathbb P_1$. Thus, one
can understand $Y_\mathrm{min}$ as a $G$-manifold,
describe the combinatorial geometry
of the steps in the Mori reduction and then reconstruct the 2:1 cover
$X\to Y$.  If $G$ is either large or has sufficiently 
complicated group structure, then the combinatorial geometry
simplifies and fine classification results can be proved.
In this regard we now mention two results from \cite {Fra}.
\begin {theorem}\label {big groups}
Let $X$ be a K3-surface and $G$ be a finite group of symplectic
automorphisms of $X$ which is 
centralized by an antisymplectic involution $\sigma $ with 
$R=\mathrm {Fix}(\sigma )\not=\emptyset $. Then, if   
$\vert G\vert >96$, it follows that the quotient $Y=X/\sigma$
is a $G$-minimal Del Pezzo surface 
and $R$ is a Riemann surface of general type.
\end {theorem}
Given our detailed knowledge of all of the surfaces 
$Y_{\mathrm {min}}$ and their automorphism groups, it would
in principle be possible to explicitly determine the K3-surfaces
which arise in this theorem. 

Although not all Mukai groups
are large in the sense of this theorem, those which are
not have a structure which is sufficiently complicated to allow 
for a precise classification.  This result can be formulated
as follows.
\begin {theorem}\label{classiMukai}
The K3-surfaces which are equipped with an
effective and symplectic action of a Mukai group $G$ centralized by an 
antisymplectic involution $\sigma $ with $\mathrm {Fix}(\sigma )\not=\emptyset $ are classified
up to equivariant equivalence in the table below.
\end {theorem}
\renewcommand{\baselinestretch}{1.25}
\begin{table}[h]
\centering
\begin{tabular}{l|l|l|l}
  & $G$ & $|G|$ & \textbf{K3-surface} $X$ \\ \hline 
1a & $L_2(7)$ & 168 & $\{x_1^3x_2+x_2^3x_3+x_3^3x_1+x_4^4 =0\} \subset \mathbb P_3$\\  \hline 
1b & $L_2(7)$ & 168 & Double cover of $\mathbb P_2$ branched along \\
& & &           $\{x_1^5x_2+x_3^5x_1+x_2^5x_3-5x_1^2 x_2^2 x_3^2 =0\}$\\ \hline
2 & $A_6$ & 360 & Double cover of $\mathbb P_2$ branched along \\
& & &           $\{10 x_1^3x_2^3+ 9 x_1^5x_3 + 9 x_2^3x_3^3-45 x_1^2 x_2^2 x_3^2-135 x_1 x_2 x_3^4 + 27 x_3^6  =0\}$\\ \hline
3a & $S_5$ & 120 & $\{\sum_{i=1}^5 x_i = \sum_{i=1}^6 x_1^2 = \sum_{i=1}^5 x_i^3=0\} \subset \mathbb P_5$\\ \hline
3b & $S_5$ & 120 & Double cover of $\mathbb P_2$ branched along $\{ F_{S_5} =0\}$\\ \hline
9 & $N_{72}$ & 72 & $\{ x_1^3+ x_2 ^3 + x_3^3 +x_4^3= x_1x_2 + x_3x_4+ x_5^2 = 0 \} \subset \mathbb P_4$ \\ \hline
10 & $M_9$ & 72 & Double cover of $\mathbb P_2$ branched along \\
& & &           $\{x_1^6+x_2^6 +x_3^6 -10(x_1^3x_2^3 + x_2^3x_3^3 +x_3^3x_1^3) =0\}$\\ \hline
11a & $ T_{48}$ & 48 & Double cover of $\mathbb P_2$ branched along $\{x_1x_2(x_1^4-x_2^4)+ x_3^6 =0\}$\\ \hline
11b & $ T_{48}$ & 48 & Double cover of $\{  x_0x_1(x_0^4-x_1^4)+ x_2^3+x_3^2=0 \} \subset \mathbb P(1,1,2,3)$\\
& &             & branched along $\{x_2=0\} $
\end{tabular}
\end{table}
\renewcommand{\baselinestretch}{1.1}
Examples 1a, 3a, 9, 10, and 11a appear in \cite{Mu} whereas the remaining provide additional examples of K3-surfaces with maximal symplectic symmetry.

There are several points concerning the above table which
need to be clarified. First, the polynomial $F_{S_5}$ in Example 3b  
can be written as
\begin {align*}
F_{S_5}=&
2(x^4yz+xy^4z+xyz^4)-2(x^4y^2+x^4z^2+x^2y^4+
x^2z^4+y^4z^2+y^2z^4)\\
&+2(x^3y^3+x^3z^3+y^3z^3)+x^3y^2z+x^3yz^2+
x^2yz^3+xy^2z^3-6x^2y^2z^2\,.
\end {align*}
The curve 
$C:=\{F_{S_5}=0\} \subset \mathbb P_2$ is singular at the points 
$[1:0:0]$, $[0:1:0]$, $[0:0:1]$ and $[1:1:1]$. The proper transform $\hat C$ of $C$ inside the Del Pezzo surface $X_5$ obtained by blowing up these four points is the normalization of $C$ and defined by a section of $-2\mathcal K_{X_5}$. The double cover $X$ of $X_5$ branched along $\hat C$ is the minimal desingularization of the double cover of $\mathbb P_2$ branched along $C$ and $X$ is a K3-surface with an action of $S_5 \times C_2$.
As is implied by Theorem \ref{big groups}
the Mori reduction with respect to the full group $G=S_5$ is such that
$X_5 = Y=Y_{\mathrm {min}}$. The
map $Y\to \mathbb P_2$ is the equivariant Mori reduction
of the Del Pezzo surface $Y$ with respect to the subgroup
$S_4$ which acts as the permutation group of the four
points which are blown up.

From the defining equation of Example 1a
one can see that this K3-surface is a $C_4$-cover of 
$\mathbb P_2$ which is branched over Klein's curve
$C$.  The preimage $\hat C$ of $C$ in the K3-surface is the
fixed point set of $C_4$. In this case $\sigma $ generates
the unique copy of $C_2$ in $C_4$ and the quotient
$X/\sigma =Y=X_2$ is a Del Pezzo surface and minimal with respect to the action of $L_2(7)$. The group
$C_4/C_2$ acts on $X_2$ and realizes $X_2$ as 2:1 cover
of $\mathbb P_2$ branched over $C$. Here $X_2$ can also be realized as the blow up $b: X_2 \to \mathbb P_2$ at the seven singular points of the sextic $\{3x^2y^2z^2 - (x^5y+y^5z+z^5x)=0\}$. Its proper transform in $X_2$ coincides with the branch locus of the map $X \to Y= X_2$. The map $b$ is the Mori reduction of $X_2$ with respect to a 
maximal subgroup $C_3\ltimes C_7$ of $L_2(7)$.

The K3-surface in Example 3a is equivariantly
embedded in $\mathbb P_5$ where $S_5$ acts by permuting the first five variables of $\mathbb C^6$ and by the character $\mathrm{sgn}$ on the sixth. 
The antisymplectic involution acts by 
$\sigma [x_1:\ldots :x_6]=[x_1:\ldots: x_5:-x_6]$.  The
quotient $X\to X/\sigma =Y=Y_{\mathrm {min}}$ is defined
by the projection $[x_1:\ldots :x_6]\mapsto [x_1:\ldots :x_5]$.
Thus $Y$ is the Del Pezzo surface of degree three defined
by the equations $\sum _1^5y_i=\sum _1^5y_i^3=0$ in
$\mathbb P_4$ known as
\emph{Clebsch cubic}. By a similar construction, Example 9 is seen to be a double cover of the \emph{Fermat cubic}. 

Finally, Example 2, the $A_6$-covering
of $\mathbb P_2$, deserves special mention.
In this
case the action of $A_6$ on $\mathbb P_2$ is given by its unique central
extension by $C_3$, which is its preimage in 
$\mathrm {SL}_3(\mathbb C)$. This was constructed by
Valentiner in the 19th century and remains of interest today
(see e.g.\,\cite {Crass}).
\subsection {Combinatorial geometry}\label{combinatorial}
The simple nature of the classification results outlined
above is at least in part due to the fact that the
possibilities for the combinatorial geometry of a Mori reduction
$Y\to Y_1\to \ldots \to Y_{\mathrm {min}}$ 
can be described in explicit ways.  An indication of this
can be found in the example in the next section.  Here we
close this section by listing the key facts which play a
role in handling this combinatorial geometry.
\begin {itemize}
\item[-]
A basic result of Nikulin (\cite {NikulinFix}) shows that the 
branch set $B$ of the covering $X \to X/\sigma =Y$,
which is the image in $Y$ of the fixed point set of 
an antisymplectic involution 
$\sigma $, is either empty, consists of two disjoint linearly equivalent elliptic curves,
is a union of rational curves or is the union of rational
curves and a Riemann surface of genus at least one.  In the
last case it is possible, and quite often happens, that
$B$ consists of only a Riemann surface of genus at least one
and no rational curves.
\item[-]
By a result due to Zhang (\cite {ZhangInvolutions}) the number of connected components of $B$
is at most ten.
\item[-]
We refer to a rational curve $E$ in $Y$ as being a \emph{Mori fiber} if it
is blown down to a point at some stage $Y_k\to Y_{k+1}$
of the reduction to a  minimal model. It can be shown that
every Mori fiber intersects the branch set in at most two
points.
\item[-]
If a Mori-fiber $E$ intersects $B$ in two points, then both
points of interestion are transversal, i.e., $E \cdot B=2$.
\item[-]
If $(E\cdot B)_p=2$, then $E\cap B=\{p\}$.
\end {itemize}
%
%
%
%
\section {The alternating group of degree six}\label {a6section}
In the previous section we considered K3-surfaces with a symplectic action of a finite group $G_\mathrm{sym}$ centralized by an antisymplectic involution, i.e., all groups under consideration were of the form $G = G_\mathrm{sym} \times C_2$. 

In this section we wish to discuss more general finite automorphims groups $\tilde G$: if $\tilde G $ contains an antisymplectic involution $\sigma$ with fixed points, then as before, we consider the quotient by $\sigma$. However, if $\sigma$ does not centralize the group $\tilde G_\mathrm{sym}$ inside $\tilde G$, the action of $\tilde G_\mathrm{sym}$ does \emph{not} descend to the quotient surface. We therefore restrict our consideration to the centralizer $Z_{\tilde G}(\sigma)$ of $\sigma$ inside $\tilde G$ and study its action on the quotient surface. 
If we are able to describe the family of K3-surfaces with $Z_{\tilde G} (\sigma)$-symmetry, it remains to identify the surfaces with $\tilde G$-symmetry inside this family. 

We consider a situation where the group $\tilde G$ contains the alternating group of degree six.
Although, a precise classification cannot be obtained at present, we achieve an improved understanding of the equivariant geometry of K3-surfaces with $\tilde G$-symmetry and classify families of K3-surfaces with $Z_{\tilde G}(\sigma)$-symmetry (cf.\,Theorem \ref{classiA6}). In this sense, this section, which
is an abbreviated version of Chapter 7 of \cite {Fra}, serves as an 
outlook on how the method of equivariant Mori reduction allows 
generalization to more advanced classification problems.
\subsection{The group $\tilde A_6$}
We let $\tilde G$ be any finite group containing the alternating group of degree six
and in the following consider a K3-surface $X$ with an effective action of $\tilde G$. 
This particular situation is considered by Keum, Oguiso, and Zhang in \cite{KOZLeech} and \cite{KOZExten} with special emphasis on the maximal possible choice of $\tilde G$: they consider a group $\tilde G = \tilde A_6$ characterized by the exact sequence
\begin{equation}\label{tilde A6}
\{\mathrm{id}\} \to A_6 \to \tilde A_6 \overset {\alpha}{\to} C_4 \to \{\mathrm{id}\}.
\end{equation}
It follows from the fact that $A_6$ is simple and a Mukai group that the group of
symplectic automorphisms $\tilde{G}_{\text{sym}}$ in $\tilde G$
coincides with $A_6$. 
Let $N := \mathrm{Inn}(\tilde{A_6}) \subset
\mathrm{Aut}(A_6)$ denote the group of inner automorphisms of $\tilde A_6$ and let $\mathrm{int} : \tilde A_6 \to N$ be the homomorphisms mapping an element $g \in \tilde A_6$ to conjugation with $g$.
 It can be shown that the group $\tilde{A_6}$ is a semidirect product $A_6
\rtimes C_4$ embedded in $N\times C_4$ by the map $(\mathrm{int}, \alpha)$
(Theorem 2.3 in \cite{KOZExten}). By Theorem 4.1 in \cite{KOZExten} the group $N$ is isomorphic to $M_{10}$
  and the isomorphism class of $\tilde A_6$ is uniquely determined by \eqref{tilde A6} and the condition that it acts on a K3-surface.

In \cite{KOZLeech} a lattice-theoretic proof of the following classification result (Theorem 5.1, Theorem 3.1, Proposition 3.5) is given. 
\begin{theorem}[\cite{KOZLeech}]
A K3 surface $X$ with an effective action of $\tilde A_6$ is isomorphic to the minimal desingularization of the surface in $\mathbb P_1 \times \mathbb P_2$ given by 
\[
S^2(X^3+Y^3 + Z^3) -3 (S^2 + T^2) XYZ =0.
\]
\end{theorem}
 The existence of an isomorphism from a K3-surface with $\tilde A_6$-symmetry to the surface defined by the equation above follows abstractly since both surfaces are shown to have the same transcendental lattice and the action of $\tilde A_6$ on the later is hidden. It is therefore desirable to obtain an explicit realization of $X$ where the action of $\tilde A_6$ is visible.

We let the generator of the
factor $C_4$ in $\tilde{A_6} = A_6
\rtimes C_4$ be denoted by $\tau$. It is nonsymplectic and has fixed points, the antisymplectic
involution $\sigma := \tau^2$ fulfils $\mathrm{Fix }_X(\sigma) \neq \emptyset$.
Since $\sigma$ is mapped to the trivial automorphism in
$\mathrm{Out}(A_6) = \mathrm{Aut}(A_6)/\mathrm{int}(A_6) \cong C_2 \times C_2$
there exists $h \in A_6$ such that $
\mathrm{int}(h) = \mathrm{int}(\sigma) \in \mathrm{Aut}(A_6)$.
The antisymplectic involution $h \sigma$ centralizes $A_6$ in $\tilde{A_6}$. 
\begin{remark}\label{valentinerremark}
If $\mathrm{Fix}_X(h \sigma) \neq \emptyset$, the K3-surface $X$ is an $A_6$-equivariant double
cover of $\mathbb{P}_2$ where $A_6$ acts as Valentiner's group and the branch locus is given by $F_{A_6}(x_1,x_2,x_3) = 10 x_1^3x_2^3+ 9 x_1^5x_3 + 9 x_2^3x_3^3-45 x_1^2 x_2^2 x_3^2-135 x_1 x_2 x_3^4 + 27 x_3^6$ (cf.\,Theorem \ref{classiMukai}). By construction, there is an evident action of $A_6 \times C_2$ on this \emph{Valentiner surface}, it is however not clear whether this surface admits the larger symmetry group $\tilde A_6$.
\end{remark}
In the following we assume that $h \sigma$ acts without fixed points on $X$ as otherwise the remark above yields an $A_6$-equivariant classification of $X$.
\subsection* {The centralizer $G$ of $\sigma$ in $\tilde{A_6}$}\label{centralizer}
We study the quotient $ \pi : X \to X/\sigma = Y$. As mentioned above, the action of the centralizer of $\sigma$ descends to an action on $Y$. We therefore start by identifying the centralizer $G :=Z_{\tilde{A_6}}(\sigma)$ of $\sigma$ in
$\tilde{A_6}$. It follows from direct computation\footnote{For this and
details of other technical arguments which have been omitted here 
see Chapter 7 of \cite {Fra}} and from the equality $\mathrm{int}(\sigma) = \mathrm{int}(h)$ that the group $G$ equals $Z_{A_6}(\sigma) \rtimes C_4$ and $Z_{A_6}(\sigma) = Z_{A_6}(h)$.
In the following, we wish to identify the group $Z_{A_6}(h)$.
Since $\mathrm{int}(\sigma) = \mathrm{int}(h)$ and $\sigma^2 = \mathrm{id}$, it follows that $h^2$ 
commutes with any element in $A_6$. As $Z(A_6)= \{ \mathrm{id} \}$, we see that $h$ is of order two. There is only one conjugacy class of elements of order two in $A_6$. We calculate $Z_{A_6}(h)$ for one particular choice of $h = (13)(24) \in A_6$. Let $c = (1234)(56)$ and $g = (24)(56)$. Then $c^2=h$ and both $c$ and $g$ centralize $h$. The group generated by $c$ and $g$ is seen to be a dihedral group of order eight;  $\langle g \rangle \ltimes \langle c \rangle = D_8 < Z_{A_6}(h)$. Now direct computations in $S_4 < A_6$ yield $\langle g \rangle \ltimes \langle c \rangle = D_8 = Z_{A_6}(h)$.

Using the assumption that $\sigma h$ acts freely on $X$ and by choosing the 
appropriate generator of $\langle c\rangle$ we find that
the action of $\tau $
on $Z_{A_6}(h)=D_8$ given by $\tau g \tau^{-1}= c^3 g$ and $\tau c\tau^{-1} = c^3$. 
Furthermore, note that the commutator subgroup $G'$ of $G$ equals
$G'=\langle c\rangle$.
\subsection*{The group $H = G / \langle \sigma \rangle$}
We consider the quotient $Y=X/\sigma $ equipped with
the action of 
$G/\sigma =:H=Z_{\tilde A_6}(\sigma)/\langle \sigma \rangle=D_8\rtimes C_2$.  The
group $C_2$ is generated by $[\tau]_\sigma$.
For simplicity, we transfer the above notation from $G$ to $H$ by writing e.g.\,$ \tau$ for
$[\tau]_\sigma$. Since $\tau g\tau ^{-1}=c^3 g= g c$, it follows that 
$H'=\langle c\rangle$. 

Let $K < G$ be the cyclic group of order eight generated by $g \tau $.
We denote the image
of $K$ in $G/\sigma $ by the same symbol.
Since $[\sigma c ]_\sigma = [c]_\sigma \in K$ it contains $H'=\langle c\rangle$
and we can write
$H=\langle \tau \rangle \ltimes K=D_{16}$.  
\begin {lemma}\label{normal groups}
There is no nontrivial normal subgroup $N$ in $H$ with
$N\cap H'=\{\mathrm{id}\}$.
\end {lemma}
\begin {proof}
If such a group exists, first consider the case $N \cap K =
\{\mathrm{id}\}$. Then $N \cong C_2$ and 
$H= K \times N$ would be Abelian, a contradiction. If $N \cap K \neq 
\{\mathrm{id}\}$ then  $N \cap K = \langle (g \tau) ^k\rangle $ for some $k \in
\{1,2,4\}$. This implies $(g \tau )^4 =c^2 \in N$ and contradicts $N \cap H' = N
\cap \langle c \rangle = \emptyset$.
\end{proof}
The following observations strongly rely the assumption that $\sigma h$ acts freely on $X$.
\begin {lemma}\label{free on B}
The subgroup $H'$ acts freely on the branch set $B = \pi(\mathrm{Fix}_X(\sigma))$
in $Y$.
\end {lemma}
\begin {proof} 
If for some $b \in B$ the isotropy group $H'_b$ is nontrivial, then $c^2(b) = h(b)=b$ and
$\sigma h$ fixes the corresponding point $\tilde b\in X$.
\end {proof}
\begin {corollary}
The subgroup $H'$ acts freely on the set $\mathcal R$ 
of rational branch curves of the covering $\pi: X \to Y$. In particular, the number of rational branch curves
$n$ is a 
multiple of four.
\end {corollary}
\begin {corollary}\label{tau-fixed}
The subgroup $H'$ acts freely on the set of $\tau $-fixed
points in $Y$.
\end {corollary}
\begin {proof}
We show $\mathrm{Fix}_Y(\tau) \subset B$.
Since $\sigma = \tau ^2$ on $X$, a $\langle \tau \rangle$-orbit $\{x, \tau x,
\sigma x, \tau^3 x \}$ in $X$ gives rise to a $\tau$-fixed point $y$ in the
quotient $Y = X /  \sigma$ if and only if $\sigma x = \tau x $. Therefore,
$\tau$-fixed points in $Y$ correspond to $\tau$-fixed points in $X$. By
definition $\mathrm{Fix}_X(\tau) \subset \mathrm{Fix}_X(\sigma)$ and the claim
follows. 
\end {proof}
\subsection {$H$-minimal models of $Y$}\label {reduction to del Pezzo}
Since $\mathrm{Fix}_X(\sigma) \neq \emptyset$, the quotient surface $Y$ is a smooth rational $H$-surface to which we apply the equivariant minimal model program. We denote by $Y_\mathrm{min}$ an $H$-minimal model of $Y$. It is known that $Y_\mathrm{min}$ is either a Del Pezzo surface or an $H$-equivariant conic bundle over $\mathbb P_1$.
\begin {theorem}\label {no equivariant fibration}
An $H$-minimal model $Y_{\mathrm {min}}$ does not admit an $H$-equivariant $\mathbb P_1$-fibration. In particular, $Y_{\mathrm {min}}$ is a Del Pezzo surface.
\end {theorem}
In order to prove this statement we begin with the following general fact which follows from the observation 
that the action of a cyclic group on a Mori fiber has two fixed points contracting to a single fixed point. 
\begin {lemma}\label{no increase}
If $Y\to Y_{\mathrm {min}}$ is an $H$-equivariant Mori reduction and $A$ a cyclic subgroup of $H$, then 
$
\vert \mathrm {Fix}_Y(A)\vert \geq 
\vert \mathrm {Fix}_{Y_{\mathrm {min}}}(A)\vert$.
\end {lemma} 
Suppose that some $Y_\mathrm{min}$ is an $H$-equivariant conic bundle, i.e., there is an $H$-equivariant fibration 
$p :Y_{\mathrm {min}}\to \mathbb P_1$ with generic fiber $\mathbb P _1$ and
let
$p _*:H\to \mathrm {Aut}(\mathbb P_1)$ 
be the associated homomorphism.
\begin {lemma}\label{ker p*}
$\mathrm {Ker}(p_*)\cap H'=\{\mathrm {id}\}$.
\end {lemma}
\begin {proof}
The elements of $\mathrm {Ker}(p_*)$ fix two points
in every generic $p $-fiber. If $h = c^2 \in H' = \langle c \rangle$ fixes
points in every generic $p$-fiber, then $h$ acts trivially on a one-dimensional subset $C \subset Y$. Since $h=c^2$ acts symplectically on $X$ it has only
isolated fixed points in $X$. Therefore, on the preimage $\tilde C = \pi^{-1}(C) \subset X$, the action of $h$ coincides with the action of $\sigma$. But then $\sigma h | _{\tilde C} = \mathrm{id}| _{\tilde C}$ contradicts the assumption that $\sigma h$ acts freely on $X$.
\end {proof}
\begin{proof}[Proof of Theorem \ref{no equivariant fibration}]
Since there are no nontrivial normal
subgroups in $H$ which have trivial intersection with 
$H'$ (Lemma \ref{normal groups}), it follows from Lemma \ref{ker p*} that $\mathrm {Ker}(p_*)=
\{ \mathrm{id} \}$, i.e., the group $H$
acts effectively on the base. 
We regard $H$ as the semidirect product 
$H=\langle \tau \rangle \ltimes K$, where $K=C_8$
is described above. 
The automorphism
$\tau$ exchanges the $K$-fixed points. We will obtain a contraction by
analyzing the $K$-actions on the fibers $F$ and $\tau F$ over its two fixed
points.  
By Lemma \ref{singular fibers of conic bundle} there are two situations which we must consider:
\begin{enumerate}
 \item 
 $F$ is a regular fiber of $Y_{\mathrm {min}}\to \mathbb P_1$.
 \item
 $F=F_1\cup F_2$ is the union of two (-1)-curves intersecting transversally in
one point. 
\end{enumerate}
We study the fixed points
of $c$, $h=c^2$ and $g \tau $ in $Y_{\mathrm {min}}$. Note that
in $X$ the symplectic transformation $c$ has precisely four fixed
points and $h$ has precisely eight fixed points. This set of eight
points is stabilized by the full centralizer of $h$, in particular by
$K$.
Since $h \sigma$ acts by assumption freely on $X$, it follows that
$\sigma$ acts freely on the set of $h$-fixed points in $X$.  If $hy=y$ for some $y \in Y$, then the preimage of $y$ in $X$ consists of two elements $x_1,\sigma x_1=x_2$. If these form an $\langle h \rangle $-orbit, then both are $ \sigma h $-fixed, a contradiction. It follows that $\{x_1,x_2 \} \subset \mathrm{Fix}_X(h)$  and the number of $h$-fixed points in $Y$ is precisely four. In particular, $h$ acts effectively on any curve in $Y$. 

Let us first consider 
case (2) where $F= F_1 \cup F_2$ is reducible.
Since $\langle c\rangle $ is the only subgroup
of index two in $K$, it follows that $\langle c \rangle $ stabilizes $F_i$ and both $c$ and $h$ have
three fixed points in $F$ (two on each irreducible component, one is the point of intersection $F_1 \cap F_2$), i.e., six fixed points on $F \cup \tau F \subset Y_\mathrm{min}$. This is contrary to Lemma \ref{no increase}
 because $h$ has at most four fixed
points in $Y_{\mathrm {min}}$.

If $F$ is regular (case (1)),
then the cyclic group $K$ has two fixed points on the rational curve
$F$. Since $h \in K$, the four $K$-fixed points on $F \cup \tau F$ are
contained in the set of $h$-fixed points on $Y_\mathrm{min}$. As
$|\mathrm{Fix}_{Y_\mathrm{min}}(h)| \leq 4$, the $K$-fixed points
coincide with the four $h$-fixed points in $Y_\mathrm{min}$, i.e., 
$\mathrm{Fix}_{Y_\mathrm{min}}(h)= \mathrm{Fix}_{Y_\mathrm{min}}(K)$.
In particular, the Mori reduction does not affect the four $h$-fixed
points $\{y_1, \dots y_4\}$ in $Y$. By equivariance of the reduction,
the group $K$ acts trivially on this set of four points. Passing to
the double cover $X$, we conclude that the action of $g \tau \in K$ on
a preimage $\{x_i , \sigma x_i\}$ of $y_i$ is either trivial or
coincides with the action of $\sigma$. In both cases it follows that
$(g \tau )^2 = c\sigma$ acts trivially on the set of $h$-fixed points in $X$. As $\mathrm{Fix}_X(c) \subset \mathrm{Fix}_X(h)$, this is contrary to the fact that $\sigma$ acts freely on $\mathrm{Fix}_X(h)$.
\end{proof}
In the following we wish to identify the Del Pezzo surface $Y_\mathrm{min}$. For this, we use the Euler characteristic
formulas,
\[
24= e(X)= 2 e(Y)-2n + \underset{\text{if $D_g$ is present}}{\underbrace{2g-2}},
\]  
where $D_g \subset B$ is a branch curve of general type, $g = g(D_g) \geq 2$, and 
$ e(Y)= e(Y_{\mathrm{min}}) + m $. Here $m = | \mathcal E|$ denotes the total number of 
Mori fibers of the reduction $Y\to Y_{\mathrm {min}}$ and $n$ denotes the 
total number of rational curves in $\mathrm{Fix}(\sigma)$.
For convenience we introduce the difference
$\delta =m -n$. 
If a branch curve $D_g$ of general type is present, then
$
13-g-\delta =e(Y_{\mathrm {min}})
$
and if it is not present
$
12-\delta =e(Y_{\mathrm {min}})
$. The inequality $e(Y_{\mathrm {min}})\le 3$ implies $m\le n+9$. 
\begin {proposition}
For every Mori fiber $E$ the orbit $H.E$ consists
of at least four Mori fibers.
\end {proposition}
\begin {proof}
A Mori fiber $E$ can either be disjoint from $B$, contained in $B$, or intersect $B$ in one or two points. We denote by $\vert HE\vert$ the number of disjoint curves in the orbit $HE$.

First assume  $E\cap B \neq \emptyset$ and $E \not\subset B$.
Since $H'$ acts freely on the branch curves
and $E$ meets $B$ in at most two points (cf.\,Section \ref{combinatorial}),
we know $\vert H'E\vert \ge 2$.
If $\vert HE\vert =2$, then the isotropy group
$H_E$ is a normal subgroup of index two which necessarily
contains the commutator group $H'$, a contradiction.

If  $E \subset B$,
we show that the $H'$-orbit of $E$ consists of four Mori fibers.
If it consisted of less than four Mori fibers, the stabilizer $H'_E \neq \{\mathrm{id}\}$ of $E$ in $H'$ would fix two points in $E \subset B$. This contradicts Lemma \ref{free on B}.

 All Mori fibers disjoint from $B$ have self-intersection (-2) and meet exactly one Mori fiber of the previous steps of the reduction in exactly one point. 
If $E \cap B = \emptyset$ there is a chain of Mori fibers $E_1, \dots, E_k =E$ connecting $E$ and $B$. 
The Mori fiber $E_1$ is the only one to have nonempty intersection with $B$ and is the 
first curve of this configuration to be blown down in the reduction process. 
The $H$-orbit of this union of Mori fibers consists of at least four copies of this chain. This is due to that fact that the $H$-orbit of $E_1$ consists of at least four Mori fibers by Case 1. In particular, the $H$-orbit of $E$ consists of at least four copies of $E$.
\end{proof}
\begin {corollary}
The difference $\delta $ is a non-negative multiple $4k$
of four. 
\end {corollary}
\begin {proof}
Above we have shown that $m$  
and $n$ are multiples of four. Therefore $\delta =4k$.
If $\delta$ was negative, i.e., 
$m < n$, there is no configuration of Mori fibers meeting the rational branch curves, which have
self-interection -4, such that the corresponding 
contractions transform them to curves on a Del Pezzo surface $Y_\mathrm{min}$ where the self-intersection
of any curve is at least -1. It follows that $\delta$ is non-negative.
\end {proof}
\begin {theorem}
Any $H$-minimal model $Y_\mathrm{min}$ of $Y$ is 
$\mathbb P_1\times \mathbb P_1$ .
\end {theorem}

\begin {proof}
If $\delta= 0$, then $n= m=0$ and $Y= Y_\mathrm{min}$. The commutator subgroup $H' \cong C_4$ acts freely on the branch locus $B$ implying $e(B)\in \{0,-8,-16, \dots \}$. Since the Euler characteristic of the Del Pezzo surface $Y$ is at least 3 and at most 11,
we only need to consider the cases $e(Y)\in \{4,8\}$.

If $\delta \neq 0$, then since $\delta \ge 4$, it follows that $e(Y_\mathrm{min})=13-g-\delta \le 7$ if a branch curve $D_g$ of general type is present, and $e(Y_\mathrm{min})=12-\delta \le 8$
if not. 

We go through the list of
of Del Pezzo surfaces with $e(Y_\mathrm{min}) \leq 8$.
\begin {itemize}
\item[-]
If $e(Y_{\mathrm {min}})=8$, then $\mathrm{deg}(Y_\mathrm{min})=4$ and $\mathrm{Aut}(Y_\mathrm{min}) = C_2^4 \rtimes \Gamma$ for $\Gamma \in \{C_2,C_4, S_3, D_{10} \}$.  If $D_{16} < C_2^4 \rtimes \Gamma$ then  $ A := D_{16} \cap C_2^4 \lhd D_{16}$ and $A$ is either trivial or isomorphic to $C_2$. In both cases $D_{16} / A$ is not a subgroup of $\Gamma$ in any of the possibilities listed above. Therefore, $e(Y) \neq 8$.
\item[-]
If $e(Y_{\mathrm {min}})=7$, then 
$\mathrm {Aut}(Y_{\mathrm {min}})=S_5$.  Since $120$ is not
divisible by $16$, we see that a Del Pezzo surface of degree five does not admit an effective action of the group $H$.
\item[-]
If $e(Y_{\mathrm {min}})=6$, then 
$A:=\mathrm {Aut}(Y_{\mathrm {min}})=
(\mathbb C^*)^2\rtimes D_{12}$.
We denote by $A^\circ \cong (\mathbb C^*)^2$ the connected component of $A$ and consider $q :A\to A/A^\circ$. Now $q (H') <  q(A)'\cong C_3$ and
$H'=C_4 < A^\circ$. We may realize $Y_{\mathrm {min}}$ as $\mathbb P_2$ blown up at the three
corner points and $A^\circ$ as the set of diagonal
matrices in $\mathrm {SL}_3(\mathbb C)$.
Every possible representation of $C_4$ in this group
has ineffectivity along one of the lines joining corner points.
But, as we have seen before, the elements of $H'$, in particular $c^2 = h$, have only isolated fixed points
in $Y_{\mathrm {min}}$.  
\item[-]
A Del Pezzo surface obtained by blowing up one or two points in $\mathbb P_2$ is never $H$-minimal and therefore does not occur
\item[-]
Finally, $Y_{\mathrm {min}}\not=\mathbb P_2$:
If $e(Y_\mathrm{min}) =3$ then either $\delta =9$ (if $D_g$ is not present), a contradiction to $\delta = 4k$, or $g + \delta =10$. In the later case, $\delta =4,8$ forces $g= 6,2$. In both cases, the Euler characteristic $2-2g$ of $D_g$ is not divisible by 4. This contradicts the fact that $H'$ acts freely on $D_g$. 
\end {itemize}
We have hereby excluded all possible Del Pezzo surfaces except $\mathbb P_1\times \mathbb P_1$ and the proposition follows. 
\end {proof}
\subsection{Branch curves and Mori fibers}
We let $M : Y \to Y_\mathrm{min} = \mathbb P_1 \times \mathbb P_1$ denote an $H$-equivariant Mori reduction of $Y$. 
\begin{lemma}
The number of Mori fibers in an $H$-orbit is at least eight.
\end{lemma}
\begin{proof}
Consider the action of $H$ on $\mathbb P_1 \times \mathbb P_1$. 
Both canonical projections are equivariant with respect to the commutator subgroup $H'= \langle c \rangle  \cong C_4$. Since $c^2 \in H'$ does not act trivially on any curve in $Y$ or $Y_\mathrm{min}$, it follows that $H'$ has precisely four fixed points in $Y_\mathrm{min} =\mathbb P_1 \times \mathbb P_1$. 
Since $h = c^2$ has precisely four fixed points in $Y$ and $\mathrm{Fix}_Y (H') = \mathrm{Fix}_Y (c) \subset  \mathrm{Fix}_Y (c^2)$, we conclude that $H'$ has precisely four fixed points in $Y$ and it follows that the Mori fibers do not pass through $H'$-fixed points. Note that the $H'$-fixed points in $Y$ coincide with the $h$-fixed points. 

Suppose there is an $H$-orbit $HE$ of Mori fibers of length strictly less then eight and let $p = M(E)$. We obtain an $H$-orbit $Hp$ in $\mathbb P_1 \times \mathbb P_1$ with $|Hp| \leq 4$. Now $| Kp| \leq 4$ implies that $K_p \neq \{\mathrm{id}\}$, in particular, $h= c^2 \in K_p$. It follows that $p$ is an $h$-fixed point. This contradicts the fact that the Mori fibers do not pass through fixed points of $h$.
\end{proof}
A total number of 24 or more Mori fibers would require 16 rational curves in $B$. This contradicts the fact that the number of connected components of the fixed point set of an antisymplectic involution on a K3-surface is at most ten (cf.\,Section \ref{combinatorial}). It therefore follows from the above lemma that the total number $m$ of Mori fibers equals 0, 8, or 16.

Recalling that the number of rational branch curves is a multiple of four, i.e., $n \in \{0,4,8\}$ and  
using the fact $m \in \{0,8,16\}$ along with $m \leq n+9$, we conclude that the surface $Y$ is of one of the following types.
\begin{enumerate}
 \item
$m=0$\\ 
The quotient surface $Y$ is $H$-minimal. The map $X \to Y \cong \mathbb P_1 \times \mathbb P_1$ is branched along a single curve $B$. This curve $B$ is a smooth $H$-invariant curve of bidegree $(4,4)$.
\item
$m=8$ and $e(Y) = 12$\\
The surface $Y$ is the blow-up of $\mathbb P_1 \times \mathbb P_1$ in an $H$-orbit consisting of eight points. 
\begin{enumerate}
 \item 
If the branch locus $B$ of $X \to Y$ contains no rational curves, then $e(B)=0$ and $B$ is either an elliptic curve or the union of two elliptic curves defining an elliptic fibration on $X$. 
\item
If the branch locus $B$ of $X \to Y$ contains rational curves, their number is exactly four (Eight or more rational branch curves of self-intersection -4 cannot be modified sufficiently and mapped to curves on a Del Pezzo surface by contracting eight Mori fibers). It follows that the branch locus is the disjoint union of an invariant curve of higher genus and four rational curves. 
\end{enumerate}
\item
$m=16$ and $e(Y) =20$\\
The map $X \to Y$ is branched along eight disjoint rational curves. 
\end{enumerate}
We may simplify the above situation by studying rational curves in $B$, their intersection with Mori fibers and their images in $\mathbb P_1 \times \mathbb P_1$.
\begin{proposition}
If $e(Y)=12$, then $n=0$.
\end{proposition}
\begin{proof}
Suppose $n \neq 0$ and let $C_i \subset Y$ be a rational branch curve. Since $C_i^2 =-4$ and $M(C_i) \subset \mathbb P_1 \times \mathbb P_1$ has self-intersection $\geq 0$ it must meet the union of Mori fibers $\bigcup E_j$. 
All possible configurations of Mori fibers yield image curves $M(C_i)$ of self-inter\-sec\-tion $\leq 4$. 
Adjunction on $\mathbb P_1 \times \mathbb P_1$
implies 
 that $g(M(C_i))=0$ and $M(C_i)$ must be nonsingular. Hence each Mori fiber meets $C_i$ in at most one point. It follows that $C_i$ meets four Mori fibers, each in one point, and $(M(C_i))^2 =0$. 
In particular, $M(C_i)$ are fibers of the canonical projections $\mathbb P_1 \times \mathbb P_1 \to \mathbb P_1$. 
The curve $C_1$ meets four Mori fibers $E_1, \dots E_4$ and each of these Mori fibers meets some $C_i$ for $i \neq 1$. After renumbering, we may assume that $E_1$ and $E_2$ meet $C_2$ and therefore $M(C_1)$ and $M(C_2)$ meet in more than one point, a contradiction. It follows that $e(Y) = 12$ implies $n=0$
\end{proof}
\begin{proposition}
If $e(Y)=20$, then $Y$ is the blow-up of $\mathbb P_1 \times \mathbb P_1$ in sixteen points $\{p_1, \dots p_{16}\} =  (\bigcup_{i=1}^4 F_i) \cap(\bigcup_{i=5}^8 F_i)$,
 where $F_1, \dots F_4$ are fibers of the canonical projection $\pi_1$ and $F_5, \dots F_8$ are fibers of $\pi_2$. The branch locus is given by the proper transform of $\bigcup F_i$ in $Y$.
\end{proposition}
\begin{proof}
We denote the eight rational branch curves by $C_1, \dots C_8$. The Mori reduction can have two steps. A slightly more involved study of possible configurations of Mori fibers shows that $0 \leq (M(C_i))^2 \leq 4$. 
As above $M(C_i)$ is seen to be nonsingular and each Mori fiber can meet $C_i$ in at most one point. Any configuration of curves with this property yields $(M(C_i))^2=0$ and $F_i = M(C_i)$ is a fiber of a canonical projection $\mathbb P_1 \times \mathbb P_1 \to \mathbb P_1$. 

If there are Mori fibers disjoint from $B$ these are blown down in the second step of the Mori reduction. Let $E_1, \dots, E_8$ denote the Mori fibers of the first step and $\tilde E_1, \dots, \tilde E_8$ those of the second step. We label them such that $\tilde E_i$ meets $E_i$. Each curve $E_i$ meets two rational branch curves $C_i$ and $C_{i+4}$ and their images $F_i = M(C_i)$ and $F_{i+4}=M(C_{i+4})$ meet with multiplicity $\geq 2$. This is contrary to the fact that they are fibers of the canonical projections. It follows that there are no Mori fibers disjoint from $B$ and all 16 Mori fibers are contrancted simultaniously. There is precisely one possible configuration of Mori fibers on $Y$ such that all rational brach curves are mapped to fibers of the canonical projections of $\mathbb P_1 \times \mathbb P_1$: The curves $C_1, \dots C_4$ are mapped to fibers of $\pi_1$ and $C_5, \dots, C_8$ are mapped to fibers of $\pi_2$. The Mori reduction contracts 16 curves to the 16 points of intersection $\{p_1, \dots p_{16} \} = (\bigcup_{i=1}^4 F_i) \cap(\bigcup_{i=5}^8 F_i)  \subset \mathbb P_1 \times \mathbb P_1$. 
\end{proof}
Let us now restrict our attention to the case where the branch locus $B$ is the union of two linearly equivalent elliptic curves and exclude this case.
\subsection*{Two elliptic branch curves}
In this paragraph we prove:
\begin{theorem}\label{two elliptic branch curves}
$\mathrm{Fix}_X(\sigma)$ is not the union of two elliptic curves.
\end{theorem}
We assume the contrary, let  $\mathrm{Fix}_X(\sigma) = D_1 \cup D_2$ with $D_i$ elliptic and let 
$f :X\to \mathbb P_1$ denote the elliptic fibration
defined by the curves $D_1$ and $D_2$.
Note that $\sigma $ acts effectively on the base $\mathbb P_1$ as otherwise $\sigma$ would act trivially in a neighbourhood of $D_i$ by a linearization argument.
It follows that the group of order
four generated by $\tau $ acts effectively on $\mathbb P_1$.

 Since the group $G$ does not contain a cyclic group of order 16, it is neither cyclic nor dihedral and therefore cannot act effectively on $\mathbb P_1$. It follows that the ineffectivity $I$ of the induced $G$-action on the base
$\mathbb P_1$ is nontrivial.
We regard
$G=C_4\ltimes D_8$ where $C_4=\langle \tau \rangle$ and $D_8$ is the centralizer of $\sigma$ in $A_6$ (cf.\,Section \ref{centralizer}) and define $J:=I\cap D_8$. 

Using explicitly the groups structure of $G$ along with the assumption that $\sigma h$ acts freely on $X$ one finds that $J$ is nontrivial.
In the following, we consider the different possibilities for the order of $J$ and show that in fact none of these occur.

If $|J|=8$ then $D_8 \subset I$. Recall that any automorphism group of an elliptic curve splits into an Abelian part acting freely and a cyclic part fixing a point. Since $D_8$ is not Abelian, any $D_8$-action on the fibers of $f$ must have points with nontrivial isotropy. This gives rise to a positive-dimensional fixed point set of some subgroup of $D_8$ on $X$, contradicting the fact that $D_8$ acts symplectically on $X$. 
It follows that the maximal possible order of $J$ is four.
\begin{lemma}\label{I does not contain c}
 The ineffectivity $I$ does not contain $\langle c\rangle$.
\end{lemma}
\begin{proof}
Assume the contrary and consider the fixed points of $c^2$. If a $c^2$-fixed point lies at a smooth point of a fiber of $f$, then the linearization of the $c^2$-action at this fixed point gives rise to a positive-dimensional fixed
point set in $X$ and yields a contradiction.
It follows that
the fixed points of $c^2$ are contained
in the singular $f $-fibers.
Since $\langle \tau \rangle $ normalizes $\langle c\rangle$ and the $\langle \tau \rangle $-orbit of a singular
fiber consists of four such fibers, we must only consider two cases:
\begin {enumerate}
\item
The eight $c^2$-fixed points are contained in four singular
fibers (one $\langle \tau \rangle $-orbit of fibers), each of these fibers contains two $c^2$-fixed points.
\item
The eight $c^2$-fixed points are contained in eight singular fibers
(two $\langle \tau \rangle$-orbits).
\end {enumerate}  
Note that $\langle c^2 \rangle$ is normal in $I$
 and therefore $I$ acts on the set of 
$\langle c^2\rangle $-fixed points. 
In the second case, all eight $c^2$-fixed points are
also $c$-fixed. This is contrary to $c$ having only four fixed
points and therefore the second case does not occur.

The first case does not occur for similar reasons: If
$c^2$ has exactly two fixed points $x_1$ and $x_2$ in some 
fiber $F$, then $\langle c\rangle $ either acts transitively
on $\{x_1,x_2\}$ or fixes both points. Since $\mathrm{Fix}_X(c) \subset \mathrm{Fix}_X(c^2)$ and $\langle c\rangle $
must have exactly one fixed point on $F$, this is impossible.
\end {proof}
\begin {corollary}
$|J| \neq 4$.
\end {corollary}
\begin {proof}
Assume $|J| = 4$. Using $\tau$ we check that no subgroup of $D_8$ isomorphic to $C_2 \times C_2$ is  normal in $G$. It follows that the group $\langle c\rangle$ is the only order four subgroup of $D_8$ which is normal in $G$ and therefore $J = \langle c\rangle$. By the lemma above this is however impossible.
\end {proof}
It remains to consider the case where $|J|=2$. 
The only normal subgroup
of order two in $D_8$ is $J=\langle h\rangle$.
\begin {lemma}
If $|J|=2$, then $I=\langle \sigma c\rangle$.
\end {lemma}
\begin {proof}
We first show that $|J|=2$ implies $|I|=4$:
If  $|I|=2$, then $ I =  \langle h \rangle$ and $G / I = C_4 \ltimes (C_2\times C_2)$. Since this group does not act effectively on $\mathbb P_1$, this is a contradiction. 
If  $|I| \geq 8$, then $G/I$ is Abelian and therefore $I$ contains the commutator subgroup $G' = \langle c \rangle$. This contradicts Lemma \ref{I does not contain c}.
It follows that $|I|=4$ and either $I \cong C_4$ or $I \cong C_2 \times C_2$. In the later case, the only possible choice is $I = \langle \sigma \rangle \times \langle h \rangle$ which contradicts the fact that $ \sigma$ acts effectively on the base. 
It follows that
$I=\langle \sigma \xi \rangle$, where $\xi ^2=h$
and therefore $\xi =c$.
\end {proof}
Let us now consider the action of $G$ on $X$ with
$I=\langle \sigma c\rangle $.
Recall that 
the cyclic group $\langle \tau \rangle$ acts effectively on the base 
and has two fixed points there. Since $\sigma =\tau ^2$, these
are precisely the two $\sigma $-fixed points. In particular,
$\langle \tau \rangle$ stabilizes both $\sigma $-fixed point 
curves $D_1$ and $D_2$ in $X$.  
Furthermore, the transformations $\sigma c$ and $c$ stabilize $D_i$ for $i =1,2$. 
Since the only fixed points of $c$ in $\mathbb P_1$ are the
images of $D_1$ and $D_2$,
$$
\mathrm {Fix}_X(c)\subset D_1\cup D_2=\mathrm {Fix}_X(\sigma).
$$
 On the other hand, we know
that $\mathrm {Fix}_X(c)\cap \mathrm {Fix}_X(\sigma )=\emptyset$. 
Thus $I=\langle \sigma c\rangle $ is not possible
and the case $|J|=2$ does not occur.

We have hereby eleminated all possibilities for $|J|$ and completed the proof of Theorem \ref{two elliptic branch curves}.

\subsection{Rough classification of $X$}
We summerize the observations of the previous section in the following classification result.
\begin{theorem}\label{roughclassiA6}
Let $X$ be a K3-surface with an effective action of the group $G$ such that $\mathrm{Fix}_X(h\sigma) = \emptyset$. Then $X$ is one of the following types:
\begin{enumerate}
 \item 
 a double cover of $\mathbb P_1 \times \mathbb P_1$ branched along a smooth $H$-invariant curve of bidegree (4,4). 
\item
 a double cover of a blow-up of $\mathbb P_1 \times \mathbb P_1$ in eight points and branched along a smooth elliptic curve $B$. The image of $B$ in $\mathbb P_1 \times \mathbb P_1$ has bidegree (4,4) and eight singular points.
\item 
 a double cover of a blow-up $Y$ of $\mathbb P_1 \times \mathbb P_1$ in sixteen points $\{p_1, \dots p_{16}\} = (\bigcup_{i=1}^4 F_i) \cap(\bigcup_{i=5}^8 F_i)$, where $F_1, \dots F_4$ are fibers of the canonical projection $\pi_1$ and $F_5, \dots F_8$ are fibers of $\pi_2$. The branch locus ist given by the proper transform of $\bigcup F_i$ in $Y$. The set $\bigcup F_i$ is an $H$-invariant reducible subvariety of bidegree (4,4).
\end{enumerate}
\end{theorem}
\begin{proof}
It remains to consider case (2) and show that the image of $B$ in $\mathbb P_1 \times \mathbb P_1$ has bidegree (4,4) and eight singular points.
We prove that each Mori fiber $E$ meets the branch locus $B$ either in two points or once with multiplicity two, i.e., we need to check that $E$ may not meet $B$ transversally in exactly one point.
If this was the case, the image $M(B)$ of the branch curve is a smooth $H$-invariant curve of bidegree $(2,2)$. 
The double cover $X'$ of $\mathbb P_1 \times \mathbb P_1$ branched along the smooth curve $M(B)= C_{(2,2)}$ is a smooth surface. Since $X$ is K3 and therefore minimal the induced birational map $X \to X'$ is an isomorphism. This is a contradiction since $X'$ is not a K3-surface.

As each Mori fiber meets $B$ with multiplicity two, the self-intersection number of $M(B)$ is 32 and $M(B)$ is a curve of bidegree (4,4) with eight singular points. These singularities are either nodes or cusps depending on the kind of intersection of $E$ and $B$. 
\end{proof}
In order to obtain a description of possible branch curves, we study the action of $H$ on $\mathbb P_1 \times \mathbb P_1$ and its invariants.
\subsection{The action of $H$ on $\mathbb P_1 \times \mathbb P_1$ and invariant curves of bidegree $(4,4)$.}
Recall that we consider the dihedral group $H \cong D_{16}$ generated by $\tau g$ of order eight and $\tau$. 
The following proposition can be obtained from direct computations:
\begin{proposition}
In appropriately chosen coordinates the action of $H$ on $\mathbb P_1\times \mathbb P_1$ is given by
\begin {align*}
c([z_0:z_1],[w_0:w_1])&=([iz_0:z_1],[-iw_0:w_1])\\
\tau ([z_0:z_1],[w_0:w_1])&=([z_1:z_0],[iw_1:w_0])\\
g([z_0:z_1],[w_0:w_1])&=([w_0:w_1],[z_0:z_1]).
\end {align*}
\end{proposition}
Given this action of $H$ on $\mathbb P_1 \times \mathbb P_1$, we wish to study the invariants and semi-invariants of bidegree $(4,4)$.
The space of $(a,b)$- bihomogeneous polynomials in $[z_0 : z_1][w_0 : w_1]$ is denoted by $\mathbb C_{(a,b)} ([z_0 : z_1][w_0 : w_1])$. 

 An invariant curve $C$ is given by a $D_{16}$-eigenvector $f \in \mathbb C_{(4,4)} ([z_0 : z_1][w_0 : w_1])$. The kernel of the $D_{16}$-representation on the line $\mathbb C f$ spanned by $f$ contains the commutator subgroup $H' = \langle c \rangle $ and $f$ is an appropiate linear combination of $c$-invariant monomials of bidegree $(4,4)$. 
It follows from the explicit form of the $H$-action that an $H$-invariant curve of bidegree $(4,4)$ in $\mathbb P_1 \times \mathbb P_1$ is 
one of the following three types
\begin{align*}
C_a &= \{a_1 f_1 + a_2 f_2 + a_3 f_3 = 0\}, \\
C_b &= \{b_1 g_1 + b_3 g_3 + b_4 g_4 =0\}, \\
C_0 &= \{g_2 =0\}.
\end{align*}
\subsection{Refining the classification of $X$}
Using the above description of invariant curves of bidegree (4,4) we may refine Theorem \ref{roughclassiA6}.
\begin{theorem}
Let $X$ be a K3-surface with an effective action of the group $G$ such that $\mathrm{Fix}_X(h\sigma) = \emptyset$. If $e(X/\sigma) = 20$, then $X/\sigma$ is equivariantly isomorphic to the blow up of $\mathbb P_1 \times \mathbb P_1$ in the singular points of the curve $C = \{f_1-f_2=0\}$ and $X \to Y$ is branched along the proper transform of $C$ in $Y$.
\end{theorem}
\begin{proof}
It follows from Theorem \ref{roughclassiA6} that $X$ is the double cover of $\mathbb P_1 \times \mathbb P_1$ blown up in sixteen points. These sixteen points are the points of intersection of eight fibers of $\mathbb P_1 \times \mathbb P_1$, four for each of fibration. 
By invariance these fibers lie over the base points $[1:1], [1:-1], [1: i], [1:-1]$ and the configurations of eight fibers is defined by the invariant polynomial $f_1-f_2$. 
The double cover $X \to Y$ is branched along the proper transform of this configuration of eight rational curves. This proper transform is a disjoint union of eight rational curves in $Y$, each with self-intersection (-4). 
\end{proof}
\begin{theorem}
Let $X$ be a K3-surface with an effective action of the group $G$ such that $\mathrm{Fix}_X(h\sigma) = \emptyset$. If $X/\sigma \cong \mathbb P_1 \times \mathbb P_1$, then after a change of coordinates the branch locus is $C_a$ for some $a_1,a_2,a_3 \in \mathbb C$.
\end{theorem}
\begin{proof}
 The surface $X$ is a double cover of $\mathbb P_1 \times \mathbb P_1$ branched along a smooth $H$-invariant curve of bidegree (4,4). The invariant (4,4)-curves $C_b$ and $C_0$ discussed above are seen to be singular at $([1:0],[1:0])$ or $([1:0],[0:1])$. 
\end{proof}
Note that the general curve $C_a$ is smooth. We obtain a 2-dimensional family $\{C_a\}$ of smooth branch curves and a corresponding family of K3-surfaces $\{X_{C_a}\}$.
It remains to consider the case (2) of the classification.
Our aim is to find an example of a K3-surface $X$ such that $X/\sigma = Y $ has a nontrivial Mori reduction $M: Y  \to \mathbb P_1 \times \mathbb P_1= Z$ contracting a single $H$-orbit of Mori fibers consisting of precisely 8 curves. In this case the branch locus $B \subset Y$ is mapped to a singular $(4,4)$-curve $C= M(B)$ in $Z$. The curve $C$ is irreducible and has precisely 8 singular points along a single $H$-orbit in $Z$. 

As we have noted above, many of the curves $C_a,C_b,C_0$ are seen to be singular at  $([1:0],[1:0])$ or $([1:0],[0:1])$. Since both points lie in $H$-orbits of length two, these curves are not candidates for our construction. This argument excludes the curves $C_b, C_0$ and $C_a$ if $a_1 = 0$ or $a_2 = 0$.

For $C_a$ with $a_3=0$ one checks that $C_a$ has singular points if and only if $a_1 = -a_2$, i.e., if $C_a$ is reducible. It therefore remains to consider curves $C_a$ where all coefficients $a_i \neq 0$. By considering the $H$-action on the irreducible component of $C_a$ one verifies that in this case $C_a$ must irreducible. We choose $a_3=1$.

One possible choice of an orbit of length eight is given by the orbit through a $\tau$-fixed point $p_\tau = ([1:1],[\pm \sqrt{i}:1])$. One checks that $p_\tau \in C_a$ for any choice of $a_i$. However, if we want $C_a$ to be singular in $p_\tau$, 
then $a_2=0$ and therefore $C_a$ is singular at points outside $H p_\tau$. It has more than eight singular points  and is therefore reducible.

All other orbits of length eight are given by orbits through $g$-fixed points $p_x= ([1:x],[1:x])$ for $x \neq 0$. One can choose coefficients $a_i(x)$ such that $C_{a(x)}$ is singular at $p_x$ if and only if $x^8 \neq 1$. If the curve $C_{a(x)}$ is irreducible, then it has precisely eight singular points $H p_x$ of multiplicity 2 (cusps or nodes) and the double cover of $\mathbb P_1 \times \mathbb P_1$ branched along $C_{a(x)}$ is a singular surface $X_\mathrm{sing}$ with precisely eight singular points. Its minimal desingularization $X$ is a K3-surface.  We obtain a diagram
\[
\begin{xymatrix}{
X_\mathrm {sing}\ar[d]^{2:1} & X \ar[d]^{2:1} \ar[l]^{\text{desing.}}\\
C_{(4,4)} 
\subset \mathbb P_1 \times \mathbb P_1 & 
\ar[l]^<<<<<{M} 
Y \supset B.}
\end{xymatrix}
\]
If $p_x$ is a node in $C_{a(x)}$, then the corresponding singularity of $X_\mathrm{sing}$ is resolved by a single blow-up. The (-2)-curve in $X$ obtained from this desingularization is a double cover of a (-1)-curve in $Y$ meeting $B$ in two points.
If $p_x$ is a cusp in $C_{a(x)}$, then the corresponding singularity of $X_\mathrm{sing}$ is resolved by two blow-ups. The union of the two intersecting (-2)-curves in $X$ obtained from this desingularization  is a double cover of a (-1)-curve in $Y$ tangent to $B$ in one point.
The information determining whether $p_x$ is a cusp or a node is encoded in the rank of the Hessian of the equation of $C_{a(x)}$ at $p_x$. The condition that this rank equals one gives a nontrivial polynomial condition. For a general irreducible member of the family $\{C_{a(x)} \, | \, x\neq 0, \,  x^8 \neq 1 \}$ the singularities of $C_{a(x)}$ are nodes.

We let $q$ be the polynomial in $x$ that vanishes if and only if the rank of the Hessian of $C_{a(x)}$ at $p_x$ is one. It has degree 24, but 16 of its solutions give rise to reducible curves $C_{a(x)}$. The remaining eight solution give rise to four different irreducible curves. These are identified by the action of the normalizer of $H$ in $\mathrm{Aut}(\mathbb P_1 \times \mathbb P_1)$ and therefore define equivalent K3-surfaces.

We summarize the discussion in the following main classification theorem.
\begin{samepage}
\begin{theorem}\label{classiA6}
Let $X$ be a K3-surface with an effective action of the group $G$ such that $\mathrm{Fix}_X(h \sigma) = \emptyset$. Then $X$ is an element of one the following families of K3-surfaces:
\begin{enumerate}
\item 
the two-dimensional family  $\{X_{C_a}\}$ for $C_a$ smooth,
\item
the one-dimensional family of minimal desingularizations of double covers of $\mathbb P_1 \times \mathbb P_1$ branched along curves in $\{C_{a(x)} \, | \, x\neq 0, \,  x^8 \neq 1 \}$. The general curve $C_{a(x)}$ has precisely eight nodes along an $H$-orbit. Up to natural equivalence there is a unique curve $C_{a(x)}$ with eight cusps along an $H$-orbit.
\item
the trivial family consisting only of the minimal desingularization of the double cover of $\mathbb P_1 \times \mathbb P_1$ branched along the curve $C_a = \{f_1-f_2=0\}$ where $a_1 =1, a_2 =-1, a_3=0$.
\end{enumerate}
\end{theorem}
\end{samepage}
\begin{corollary}
 Let $X$ be a K3-surface with an effective action of the group $\tilde A_6$. If $\mathrm{Fix}_X(h \sigma) = \emptyset$, then $X$ is an element of one the families (1) -(3) above. If $\mathrm{Fix}_X(h \sigma) \neq \emptyset$, then $X$ is $A_6$-equivariantly isomorphic to the Valentiner surface.
\end{corollary}
\subsection{Summary and outlook}
Our initial goal in this section was the description of K3-sur\-fa\-ces with $\tilde A_6$-symmetry. Using the group structure of $\tilde A_6$ this problem is now  divided into two possible cases corresponding to the question whether $\mathrm{Fix}_X(h\sigma)$ is empty or not. If it is nonempty, the K3-surface with $\tilde A_6$-symmetry is the Valentiner surface (Remark \ref{valentinerremark}).
If is is empty, our discussion in the previous sections has reduced the problem to finding the $\tilde A_6$-surface in the families of surfaces $X_{C_a}$ with $D_{16}$-symmetry.
It is known that a K3-surface with $\tilde A_6$-symmetry has maximal Picard rank 20. This follows from a criterion due to Mukai (\cite{Mu}) and is explicitely shown in \cite{KOZLeech}. 
All surfaces $X_{C_a}$ for $C_a \subset \mathbb P_1 \times \mathbb P_1$ a (4,4)-curve are elliptic since the natural fibration of $\mathbb P_1 \times \mathbb P_1$ induces an elliptic fibration on the double cover (or is desingularization). 

A possible approach for finding the $\tilde A_6$-example inside our families is to find those surfaces with maximal Picard number by studying the elliptic fibration.
It would be desirable to apply the following formula for the Picard rank of an elliptic surface $f: X \to \mathbb P_1$ with a section (cf.\,\cite{shiodainose}):
\[
 \rho(X) = 2 + \mathrm{rank}(MW_f) + \sum_i (m_i-1),
\]
where the sum is taken over all singular fibers, $m_i$ denotes the number of irreducible components of the singular fiber and $\mathrm{rank}(MW_f)$ is the rank of the Mordell-Weil group of sections of $f$. 

First, one has to ensure that the fibration under consideration has a section. One approach to find sections is to consider the quotient $q:\mathbb P_1 \times \mathbb P_1 \to \mathbb P_2$ and the image of the curve $C_a$ inside $\mathbb P_2$. For an appropiate bitangent to $q(C_a)$ its preimage in the double cover of $\mathbb P_1 \times \mathbb P_1$ is reducible and both its components define sections of the elliptic fibration. For the special curve $C_a$ with eight nodes the existence of a section (two sections) follows from an application of the Pl\"ucker formula to the curve $q(C_a)$ with 3 cusps and its dual curve. 

As a next step, one wishes to understand the singular fibers of the elliptic fibrations. Singular fibers occur whenever the branch curve $C_a$ intersects a fiber $F$ of the $\mathbb P_1 \times \mathbb P_1$ in less than four points. Depending on the nature of intersection $F \cap C_a$ one can describe the corresponding singular fiber of the elliptic fibration. For $C_a$ the curve with eight cusps one finds precisely eight singular fibers of type $I_3$, i.e., three rational curves forming a closed cycle. In particular, the contribution of all singular fibers $\sum_i (m_i-1)$ in the formula above is 16. In the case where $C_a$ is smooth or has eight nodes, this contribution is less. 

In order to determine the number $ \rho(X_{C_a})$ it is neccesary to either understand the Mordell-Weil group or to find curves giving additional contribution to $\mathrm{Pic}(X_{C_a})$ not included in $2 + \sum_i (m_i-1)$. 

In conclusion, the method of equivariant Mori reduction applied to quotients $X/\sigma$ yields an explicit description of families of K3-surfaces with $D_{16} \times \langle \sigma \rangle$-symmetry and by construction, the K3-surface with $\tilde A_6$-symmetry is contained in one of these families.  It remains to find criteria to characterize this particular surface 
inside these families. The possible approach by understanding the function  $a \mapsto \rho( X_{C_a})$
using the elliptic structure of $X_{C_a}$ requires a detailed analysis of the Mordell-Weil group. 
%
%

%
\end {document}